\documentclass[12pt]{article}
\usepackage{amsmath}
\usepackage{graphicx}
\usepackage{enumerate}
\usepackage{natbib}
\usepackage{url} % not crucial - just used below for the URL 
\usepackage{bm}
\usepackage{amsfonts}
\usepackage{mathtools}
\usepackage{amsthm}

\newtheorem{theorem}{Theorem}[section]
\newtheorem{corollary}{Corollary}[theorem]

\newtheorem{definition}{Definition}[section]
\usepackage{xcolor} 
\usepackage{mathrsfs}
\newcommand{\blind}{1}

\begin{document}

\def\spacingset#1{\renewcommand{\baselinestretch}%
	{#1}\small\normalsize} \spacingset{1}

%%%%%%%%%%%%%%%%%%%%%%%%%%%%%%%%%%%%%%%%%%%%%%%%%%%%%%%%%%%%%%%%%%%%%%%%%%%%%%

\if1\blind
{
	\title{\bf Graph-Based Equilibrium Metrics for Dynamic Supply-Demand Systems with Applications to Ride-sourcing Platforms  \footnote{
This work was done when Dr. Hongtu Zhu took the leave of absence from the University of North Carolina at Chapel Hill. 
 The readers are
welcome to request reprints from  Dr. Hongtu Zhu. Email: bowenhongtu@gmail.com; Phone: 919-966-7272.
 }}
	\author{Fan Zhou$^1$, 		Shikai Luo$^2$,  Xiaohu Qie$^2$,  Jieping Ye$^2$,  and Hongtu Zhu$^2$ \\%\thanks{
		%The authors gratefully acknowledge \textit{please remember to list all relevant funding sources in the unblinded version}}
 $^1$Shanghai University of Finance and Economics
	and		 $^2$Didi Chuxing }
\date{} 
	\maketitle
} \fi

\if0\blind
{
 
	\begin{center}
		{\LARGE\bf Graph-Based Equilibrium Metrics for Dynamic Supply-Demand Systems}
	\end{center}
	\medskip
} \fi

\vskip -0.50 in 
\begin{abstract} 
How to dynamically measure the local-to-global spatio-temporal coherence between demand and supply networks is a fundamental task for ride-sourcing platforms, such as DiDi. Such coherence measurement is critically important for the quantification of the market efficiency and the comparison of different platform policies, such as dispatching. The aim of this paper is to introduce a graph-based equilibrium metric (GEM) to quantify the distance between demand and supply networks based on a weighted graph structure. We formulate GEM as the optimal objective value of an unbalanced transport problem, which can be efficiently solved by optimizing an equivalent linear programming.  We examine how the GEM can help solve three operational tasks of ride-sourcing platforms. The first one is that GEM achieves up to 70.6$\%$ reduction in root-mean-square error  over the second-best distance measurement for the prediction accuracy. The second one is that the use of GEM for designing order dispatching policy increases  answer rate and drivers' revenue for more than 1$\%$,   representing a huge improvement in number. The third one is that GEM is to serve as an endpoint for comparing different platform policies in AB test.
\end{abstract}
\noindent%
\vfill

\newpage
\spacingset{1.5} % DON'T change the spacing!
\section{Introduction}
\label{sec:intro}
\begin{figure}[h]
\vskip -0.3 in 
	\centering
	\includegraphics[width=0.5\linewidth]{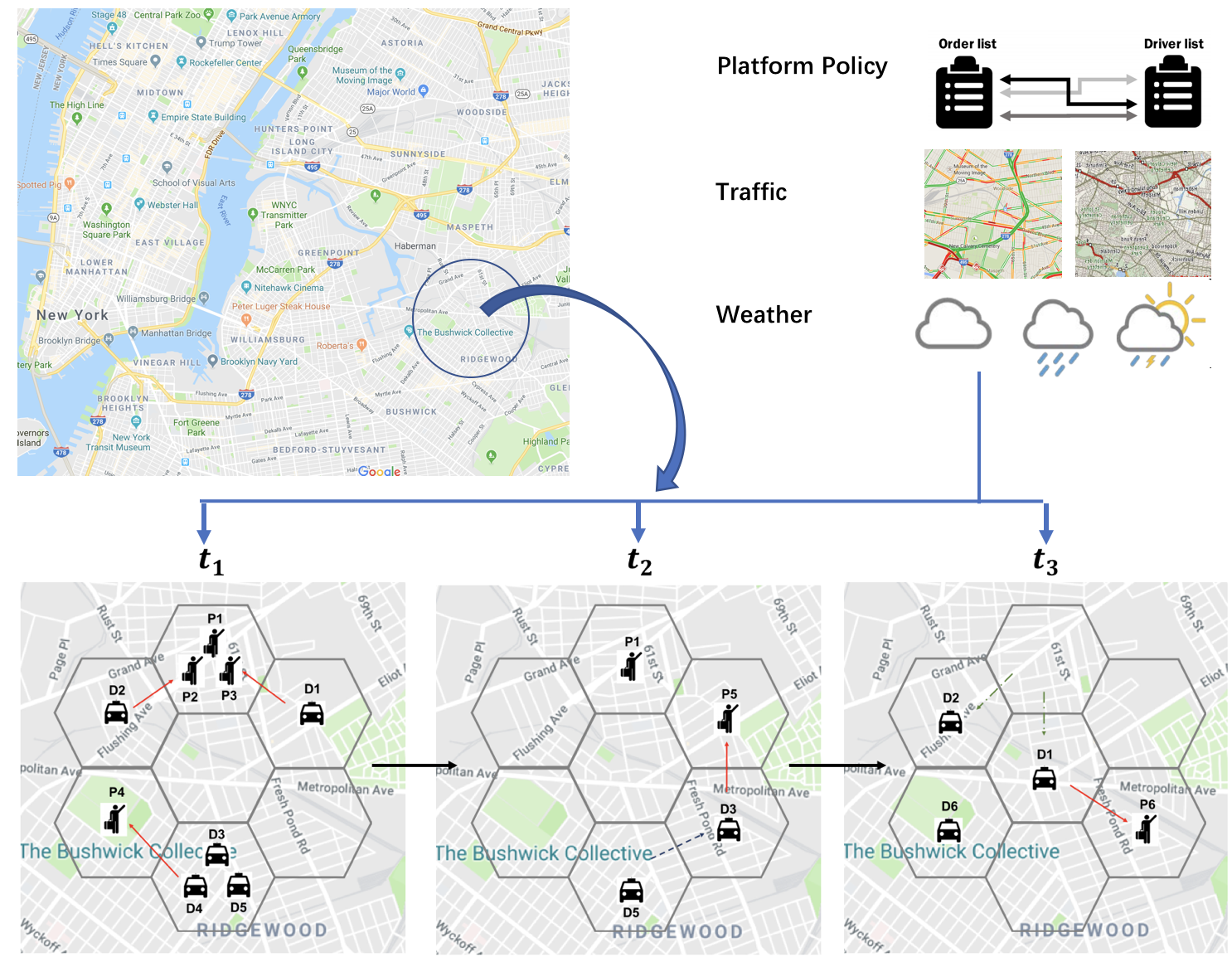}
	\caption{Dynamic supply and  demand networks at three time points in a representative ride-sourcing platform. We  divide the whole city into multiple 
		hexagon areas.  }
	\label{fig:big_picture} 
\end{figure} 

Large volumes of data collected from multiple spatio-temporal  networks are increasingly studied in diverse fields including climate science, social sciences, neuroscience, epidemiology, and transportation. In addition, those spatio-temporal networks may interact with each other across spatial and/or temporal dimension. A typical example is  that the dynamic  demand and supply networks of a ride-sourcing platform \citep{WangYang2019} are two sequences of un-normalized masses measured on the same undirected (or directed) graph $G=(\mathbb V, \mathbb E)$, where $\mathbb V$ and $\mathbb E$ are, respectively, a vertex set and  
a set of edges connecting vertex pairs. Figure \ref{fig:big_picture} illustrates how the two complicated networks interact with each other and evolve over  time. Specifically, 
a  city  is  divided into hundreds of non-overlapping grids as the vertex set $\mathbb V$ with the edge structure $\mathbb E$ determined by road networks and location functionalities. Both demands and supplies are observed across grids at each time window with possibly  different total masses and  distributions. The ride-sourcing platform uses some order  dispatching 
policy to   match  customer requests with possible surrounding idle drivers, while  after finishing serving  assigned orders, drivers return back to the supply pool to prepare for the next feasible matching. The aim of this paper is to 
address a  fundamental   question of interest for the demand and supply networks of   two-sided markets.

%A fundamental question  of interest for ride-sourcing platforms is 
%\begin{itemize} 
% \item{} how to evaluate whether a ride-sourcing platform achieves  a healthy equilibrium between dynamic supply and demand networks, what we call market equilibrium.  
% \end{itemize} 
% A healthy market equilibrium needs to meet three basic conditions.  The first is that  each  passenger's request  could be easily and quickly satisfied.  Second, the total  idle time for all drivers could be minimal. The last is to minimize various `transporting' costs, such as wait time for riders and pick-up time for drivers, which depend on the size and city geography, weather conditions and market competition. Achieving  such healthy market equilibrium  requires the precise  operation of various platform policies including dispatching, dispositioning, pricing, and subsiding under a given non-platform environment, including weather, holiday, economic level,  and government policy.   
%   how to achieve this equilibrium.  

The  fundamental question  of interest   that we consider here  is how  to quantify the spatial equilibrium of dynamic supply-demand networks for  
 two-sided markets, particularly ride-sourcing platforms (e.g., Uber and DiDi). 
To solve this question, we first  introduce a  weighted graph structure $(G,W,C)$ to characterize the 
transport  network and transport  costs of  a city.    Specifically, we divide each market into $N$ 
disjoint areas and regard them as vertices, denoted as $\mathbb V=\{v_1, \ldots, v_N\}$. 
Let  $\mathbb E$  be a set of edges between any possible pair of vertices such that    $(v_i, v_j) \in \mathbb E \subset \mathbb V\times \mathbb V$ is an edge  equipped with an nonegative weight $w_{ij}$ (e.g., transportation cost).   For all $(v_i, v_j) \notin \mathbb E,$ we set  $w_{ij} = \infty$. 
The weighted graph structure consists of  
an undirected (or directed) graph $G=(\mathbb V, \mathbb E)$ as well as a weight matrix $W=(w_{ij})$, where $w_{ij}$s' are nonegative weights.   
A graph-based  transport  cost   from $v_i$ to $v_j$ is defined as 
 $c_{ij} = \min_{K \geq 0, (i_k)_{k=0}^K:v_i \rightarrow v_j} \{\sum_k w_{i_k, i_{k+1}}: \forall k \in [\![0, K-1]\!], (v_{i_k}, v_{i_{k+1}}) \in \mathbb E \},  
$
where $(i_k)_{k=0}^K: v_i \rightarrow v_j$ denotes any  path on $G$  through $\mathbb E$ starting from $v_{i_0} = v_i$ and ending at $v_{i_K} = v_j$. Thus,  
$c_{ij}$ is the geodesic distance from $v_i$ to $v_j$ or the minimal cost of transporting one unit of object from $v_i$ to $v_j$. 
Thus,   we can define a  transport cost matrix on $(G, W)$, denoted as  $ C=(c_{ij}) \in R^{N \times N}$. The  $C$ may be time variant, since it depends  on  the real-time traffic and weather conditions for ride-sourcing platform.  The  $C$ is possibly  asymmetric since  the graph $G$ can be directed. 
% Throughout the paper, 
%we consider the weighted graph structure  $(G, W, C)$. 

Second, we need to introduce 
a distance (or metric) to quantify the difference between  demand and supply masses   at each time interval and across time on $(G,W,C)$.   
At a given time interval, we define  $\nu_j = \nu(v_j)$ and $\mu_j = \mu(v_j)$ as the point masses at vertex $v_j$ for the two measures  $\nu$ and  $\mu$, which, respectively, represent the number of customer requests and available drivers inside the vertex $v_j$ of the ride-sourcing platform \citep{WangYang2019}. 
The supply and demand systems at each timestamp can be modeled as two discrete Lebesgue measures $\mu$ and $\nu$ on    $(G, W, C)$ with locally finite masses such that  $\max(\mu(\mathbb V_0), \nu(\mathbb V_0))$ is finite for every compact set $\mathbb V_0 \subset \mathbb V$.  We consider a general case that the two measures can be unbalanced, that is,  $\mathbf{\mu} = \sum_{i=1}^N \mu_{i}$ and  $\mathbf{\nu} = \sum_{i=1}^N \nu_{i}$ may be unequal to each other. Defining a metric between $\mu$ and $\nu$ falls into the  field of optimal transport.

Optimal transport has been widely  studied in diverse disciplines, such as statistics, applied mathematics,    medical imaging,  
    and computer vision. Wasserstein-based metrics  based on the mathematics of optimal mass transport  have been  proved to be powerful tools for comparing 
objects in complex spaces.  Some successful applications include solving transport partial differential equation (PDE) \citep{ambrosio2008hamiltonian},  imaging processing \citep{rabin2015convex}, statistical inference in machine learning \citep{solomon2014wasserstein}, manifold diffeomorphisms \citep{Grenander2007},  and serving as the cost function for training Generative Adversarial Networks  \citep{arjovsky2017wasserstein}, among many others. However, existing Wasserstein-based metrics  are not directly applicable to the comparison of two unbalanced measures  defined 
%discrete measures defined 
on $(G,W,C)$  as detailed in Section \ref{sec21}.

We introduce  a   graph-based  equilibrium metric (GEM) and formulate it 
as an unbalanced optimal transport problem.   
Our main   contributions  are summarized as follows. 
First, we propose a  novel GEM, which  can be regarded as a restricted generalized Wasserstein distance,  to quantify the distance between dynamic demand and supply networks  on the weighted graph structure.   
	It not only allows the optimal transport guided by asymmetric costs and node connections,  but also accounts for unbalanced masses. 
	It also allows one of the two sides (supplies) to play the transporting role and the other (demands) to be fixed, which satisfies the physical interpretation of ride-sourcing platforms.  
	Second, varying the size of each vertex leads to 
	multilevel GEMs and their corresponding optimal transport functions. At the finest scale, our GEM reduces to solving an unbalanced assignment problem and 
	its corresponding optimal transport function contains  many local details. 
	In contrast, at a relatively coarse scale,  it gives  a coarse representation (or low frequency patterns) of  
	the optimal transport function.   
Third, 
numerically,    the calculation of  GEM can be reformulated as a standard linear programing (LP) problem. 
	Theoretically, we  investigate several theoretical properties of GEM including the convergence of the LP algorithm for computing GEM, the expectation of GEM, 
	and the metric property, additive property and   weak convergence of GEM.  
Fourth,  we apply GEM to the `supply-demand diagnostic data set' obtained from  the DiDi Chuxing  in order  to address  some  important operational tasks,  such as   the prediction of  the 
 	efficiency of  a given dispatching policy.

The remainder of this paper is structured as follows. 
	Section   \ref{sec2} develops   graph-based  equilibrium metrics and their computational approach, while discussing  their potential applications. 
	Section \ref{sec3} studies four theoretical properties associated with  GEM.  
	Section \ref{sec4} demonstrates the  applications of GEM in the intelligent operations of DiDi Chuxing.

%Body of paper.  Margins in this document are roughly 0.75 inches all
%around, letter size paper.

%\begin{figure}
%\begin{center}
%\includegraphics[width=3in]{fig1.pdf}
%\end{center}
%\caption{Consistency comparison in fitting surrogate model in the tidal
%power example. \label{fig:first}}
%\end{figure}

%\begin{table}
%\caption{D-optimality values for design $X$ under five different scenarios.  \label{tab:tabone}}
%\begin{center}
%\begin{tabular}{rrrrr}
%one & two & three & four & five\\\hline
%1.23 & 3.45 & 5.00 & 1.21 & 3.41 \\
%1.23 & 3.45 & 5.00 & 1.21 & 3.42 \\
%1.23 & 3.45 & 5.00 & 1.21 & 3.43 \\
%\end{tabular}
%\end{center}
%\end{table}

\vskip -0.5 in 

\section{Methodologies}  \label{sec2}

\vskip-0.2 in

\subsection{Existing  Wasserstein-type  
	Distances}
\label{sec21}

\begin{figure}
	\centering
	\includegraphics[width=0.6\linewidth]{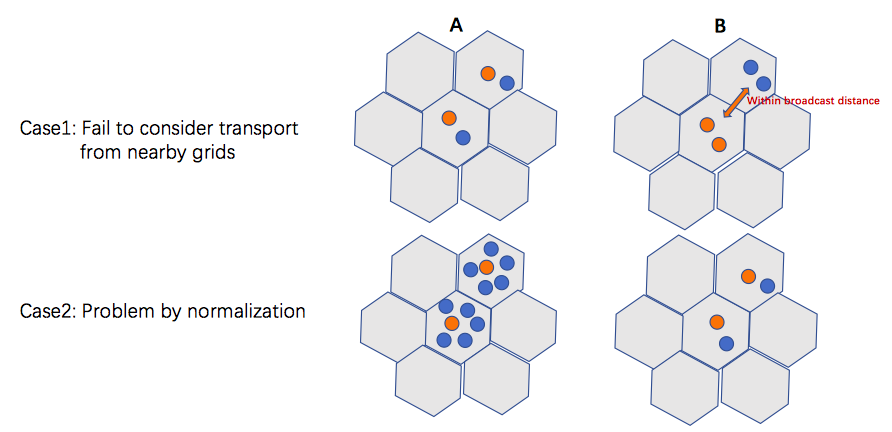}
\vskip -0.3 in  
	\caption{ Examples of supply-demand networks for illustrating the limitations of Hellinger distance:   blue and orange objects represent orders and idle drivers, respectively.
 The first row illustrates the consequence of not   transporting   objects into nearby  vertices. Specifically, the left one has a zero Hellinger distance, whereas the right one has a nonzero Hellinger distance without transporting the two drivers into the  nearby vertex with two orders. The second row illustrates that the normalization step can remove the imbalance between orders and idle drivers.  Specifically, the left one has a large imbalance between 
orders and idle drivers, whereas after the normlization,  the right one has a zero Hellinger distance.   }
	\label{fig:hellinger}
\vskip -0.2 in 
\end{figure}

Many approaches have been proposed to measure the distance between two measures (or distributions) on a metric space. Most of them fall into two broad categories including the aggregation of pixel-wise differences and the  transport cost of moving one measure to match the other. Measurements in the first category include the $L_p$-distance,  the Total Variation (TV) distance,  and the Kullback-Leibler (KL) divergence \citep{cha2007comprehensive}, among others. A  typical example is the Hellinger Distance  reviewed as follows.
\begin{definition}
	\textbf{(Hellinger Distance)} Let $M(X)$  and $M_{+}(X)$ be the vector space of Radon measures  and the cone of nonnegative Radon measures on a Hausdorff topological space $X$.  
	Then, we use $\mu$ and $\nu$ to denote two probability measures that are absolutely continuous with respect to a third probability measure $\nu_0$. The square of the Hellinger distance between $\mu$ and $\nu \in M_{+}(X)$ is defined as 
$
	D_H^{2}(\mu, \nu)=0.5  \int_X \left({\sqrt {d\mu/d\nu_0 }}-{\sqrt {d\nu/d\nu_0 }}\right)^{2}d\nu_0, 
$	where $d\mu/d\nu_0$ and $d\nu/d\nu_0$ are the Radon-Nikodym derivatives of $\mu$ and $\nu$, respectively.
\end{definition}
All these metrics suffer from two  major issues. Please refer to Figures \ref{fig:hellinger} for details.  
First, all these metrics not only fail to consider the connections 
among different locations (vertices in graph), but also ignore the topological (or geometric) structure of $X$. Second, the use of Hellinger-type distances requires a normalization  step to enforce  $\mu(X)
=\int_X d\mu=\nu(X)=\int_X d\nu$, which can create  
a false balance issue.

To address these two issues, the second category of  metrics, such as the Wasserstein distance \citep{villani2008optimal},  is proposed by solving an optimal transport problem. In the real world,  original supply resources can usually be  transported to achieve a better equilibrium  between $\mu$ and $\nu$. 
%the metric measures the minimum "cost" of turning $\mu$ into $\nu$.}
All those  distances   have deep connections to well studied assignment  
problems  from combinatorial optimization \citep{Steele1987}.

\begin{definition}
	\textbf{(Wasserstein Distance)} Let $X$ and $Y$ be  Hausdorff topological spaces and $X\times Y$ be their product space.    We introduce a lower semi-continuous  function 
	$c: X \times Y \rightarrow R\cup \{\infty\}$,  an nonnegative measure (or a transport function) $\gamma \in M_{+}(X \times Y)$,    and an equality constraint $\iota_{\{=\}}(\alpha|\beta)$ which is $0$ if $\alpha = \beta$ and $\infty$ otherwise. Then the optimal transport problem for measures $\mu \in M_{+}(X)$ and $\nu \in M_{+}(Y)$  with
	the same total masses, that is $\mu(X)=\nu(Y)$,  can be defined as 
	\begin{equation}\label{Was}
	D_W(\mu, \nu|c) = \inf_{\gamma \in M_{+}(X \times Y)} \left\{\int_{X \times Y} c d\gamma + \iota_{\{=\}}(P^X_{\#}\gamma|\mu) + \iota_{\{=\}}(P^Y_{\#}\gamma|\nu)\right\}, 
	\end{equation}
	where  
	$P^X_{\#}\gamma$ and $P^Y_{\#}\gamma$ denote the first and second marginals of $\gamma$, respectively.
\end{definition}

Intuitively, $\gamma$ denotes a transport plan, measuring how far you have to move the mass of $\mu$ to turn it into $\nu$.
Standard  optimal transport  in (\ref{Was}) is only meaningful  whenever $\mu$ and $\nu$ have the same total masses.
Whenever $\mu(X)\not=\nu(Y)$, there is no feasible $\gamma$ in (\ref{Was}). 
For the real-world ride-sourcing platforms, however, it is important to compute some sort of relaxed transportation between two arbitrary non-negative  measures. 
An improved approach is to build an unbalanced optimal transport problem by introducing   two divergences over $X$ and $Y$, denoted as $\mathcal D_{\varphi_1}$ and $\mathcal D_{\varphi_2}$,  respectively \citep{chizat2018scaling, liero2018optimal}. 
%The formal definition of the $\varphi$-divergence is given as follows. 
\begin{definition}
\textbf{(Divergences)}. Let $\varphi$ be an entropy function. For $\mu, \nu \in M(T)$, $\frac{d\mu}{d\nu}\nu + \mu^\bot$ is the Lebesgue decomposition of $\mu$ with respect to $\nu$. The divergence $D_\varphi$ is defined by
$D_\varphi(\mu|\nu) := \int_T \varphi (\frac{d\mu}{d\nu})d\nu +\varphi_{\infty}'\mu^\bot(T)$ 
if $\mu$ and $\nu$ are nonnegative and $\infty$ otherwise. 
\end{definition}
Now, we can give the formal definition of Generalized Wasserstein Distance. 
\begin{definition}
	\textbf{(Generalized Wasserstein Distance (GWD))} Let $c: X \times Y \rightarrow [0, \infty]$ be a lower semi-continuous  function, the unbalanced optimal transport problem is 
	\begin{equation}\label{geneWas}
	D_{\varphi_1, \varphi_2}(\mu, \nu|c) = \inf_{\gamma \in M_{+}(X \times Y)} \left\{\int_{X \times Y} c d\gamma + \mathcal D_{\varphi_1}(P^X_{\#}\gamma|\mu) + \mathcal D_{\varphi_2}(P^Y_{\#}\gamma|\nu)\right\}.  
	\end{equation}
\end{definition} 
Different from standard Wasserstein Distance which normalizes the input measures into probability distributions, GWD quantifies in some way the deviation of the marginals of the transport plan $\gamma$ from the two unbalanced measures $\mu$ and $\nu$
%assigned measures 
by using $\varphi$-divergence.
Although $D_{\varphi_1, \varphi_2}(\mu, \nu|c)$ enjoys some nice properties, such as metric property  \citep{chizat2018scaling, liero2018optimal},  the solution to (\ref{geneWas}), denoted as $\gamma_*$, may not have any physical meaning. For the ride-sourcing business, such $\gamma_*$ is critically important for assigning supplies to demands since it can be regarded as 
the graph representation of a dispatching policy. Therefore, the use of $D_{\varphi_1, \varphi_2}(\mu, \nu|c)$  still cannot fully cover the 'useful' relative size between   $\mu$ and $\nu$, since
it may underestimate  unmatched resources by allowing some infeasible transports, that is, the space $M_{+}(X\times Y)$ is too large to be useful.  
Three major issues of using $D_{\varphi_1, \varphi_2}(\mu, \nu|c)$  are given as follows.  
	The first issue is that in many applications (e.g., ride-sourcing platform), point masses in only one of the two measures are allowed to be transported and those in the other measure are fixed. In this case,  the symmetric property does not hold.  
	The second issue is that neither $W$ nor $C$ can be used to define a standard metric space on $G=(\mathbb V, \mathbb E)$, since transport  cost (or weight) matrix  may not satisfy  the three key assumptions of standard metrics.  
	For instance, the transport cost  from $v_i$ to $v_j$ may be unequal to that from $v_j$ to $v_i$, since transport cost matrix $C \in R^{N \times N}$ can be asymmetric  for directed graphs. Moreover,  the direct  transport cost  from $v_i$ to $v_j$  may be larger than or equal to the sum of  the  transport cost from $v_i$ to $v_k$ and that from 
	$v_k$ to $v_j$.
	The third issue is that  in some  applications, such as supply-demand networks,   the  transport cost from $v_i$ to $v_j$ may not be a constant and the  transport cost  from a vertex to itself may not be zero.
	It is possible that supply units  at vertex $v_i$ have  their individual  transport costs of moving within/outside the vertex $v_i.$ Subsequently, their transport costs  from $v_i$ to $v_j$  may follow a distribution instead of being a constant.

\vskip-0.2 in 

\subsection{Graph-based Equilibrium Metrics}  \label{sec22}

On    $(G, W, C)$, we formally introduce our GEMs for two discrete measures $\mu$ and $\nu$ in $M_{+}(\mathbb V)$, among which point masses in $\mu$  are allowed to be transported and those in $\nu$ are fixed.     We need to introduce some 
notations.   In this case, we have $X=Y=\mathbb V$ and use 
$P^V_{\# 1}\gamma$ and $P^V_{\# 2}\gamma$ to represent $P^X_{\#}\gamma$ and  $P^Y_{\#}\gamma$, respectively.  
Let $|\mu|=\sum_{i=1}^N\mu(v_i)$ and $|\mu-\tilde \mu|=\sum_{i=1}^N|\mu(v_i)-\tilde\mu(v_i)|$.  
For $i=1, \ldots, N$, we use ${\mathcal N}_i$ to denote the neighboring set   of 
$v_i$ in $\mathbb V$, which contains $v_i$ and its (possibly high-order)  neighboring vertexes. Moreover,  $v_i \in \mathcal N_j$ does not ensure $v_j \in \mathcal N_i$ since the traffic and road networks may constrain directly  transporting cars from $v_j$ to $v_i$.

Let $c: \mathbb V \times \mathbb V \rightarrow R\cup \{\infty\}$ be a   function and  $\gamma \in M_{+}(\mathbb V \times \mathbb V)$ be a nonnegative measure.  
The general form of our GEM on   $(G, W, C)$ is written as 
\begin{equation}\label{em}
\rho_{\lambda} (\mu,\nu | G, C) = \inf_{\tilde \mu \in M_{+}(\mathbb V), \gamma\in M_{+}(\mathbb V\times \mathbb V)} \left\{|\nu - \tilde \mu| + \lambda 
\int_{\mathbb V \times \mathbb V} c d\gamma  \right\}  
\end{equation}
%and the two transport constraints in (\ref{eq4}) are modified to
%\begin{equation}\label{eq8}
%(P^X_{\#}\gamma)_i(x) = \int_{N_i} \gamma \, dy \,\,\,\,\, \text{and} %\,\,\,\,\, (P^X_{\#}\gamma)_i(y) = \int_{N_i} \gamma \, dx
%\end{equation}
subject to  an equality constraint
and two sets of transport constraints given by 
%\begin{equation}\label{eq8}
%|\mu|=|\tilde\mu|,~~~~ 
%(P^V_{\# 1}\gamma)(v_i) = \sum_{v \in \mathcal N_i} \gamma(v_i, v)  = \mu_i \,\,\,\,\, \text{and} \,\,\,\,\, (P^V_{\# 2}\gamma )(v_i) = \sum_{v \in \mathcal N_i} \gamma(v, v_i) \,  = \tilde \mu_i, 
%\end{equation}
\begin{equation}\label{eq8}
|\mu|=|\tilde\mu|,~~~~ 
(P^V_{\# 1}\gamma)(v_i) = \sum_{v_j \in \mathcal N_i} \gamma(v_i, v_j)  = \mu_i \,\,\,\,\, \text{and} \,\,\,\,\, (P^V_{\# 2}\gamma )(v_i) = \sum_{v_i \in \mathcal N_j} \gamma(v_j, v_i) \,  = \tilde \mu_i, 
\end{equation}
where $\lambda$ is a non-negative hyper-parameter.   The three sets of constraints in (\ref{eq8}) ensure that $\tilde\mu$   shares the same total mass 
with $\mu$ and $\gamma$ transports  $\mu$ to $\tilde\mu$.  Thus, the feasible  set for (\ref{eq8}) is much smaller than that for (\ref{geneWas}). 
The integration of $\lambda 
\int_{\mathbb V \times \mathbb V} c d\gamma $ and the three sets of  constraints in (\ref{eq8}) is equivalent to the balanced Wasserstein distance in (\ref{Was}), so GEM is the integration of 
the  balanced Wasserstein distance and the $L_1$ norm. 

\begin{figure}[h]
	\centering
	\includegraphics[width=0.7\linewidth]{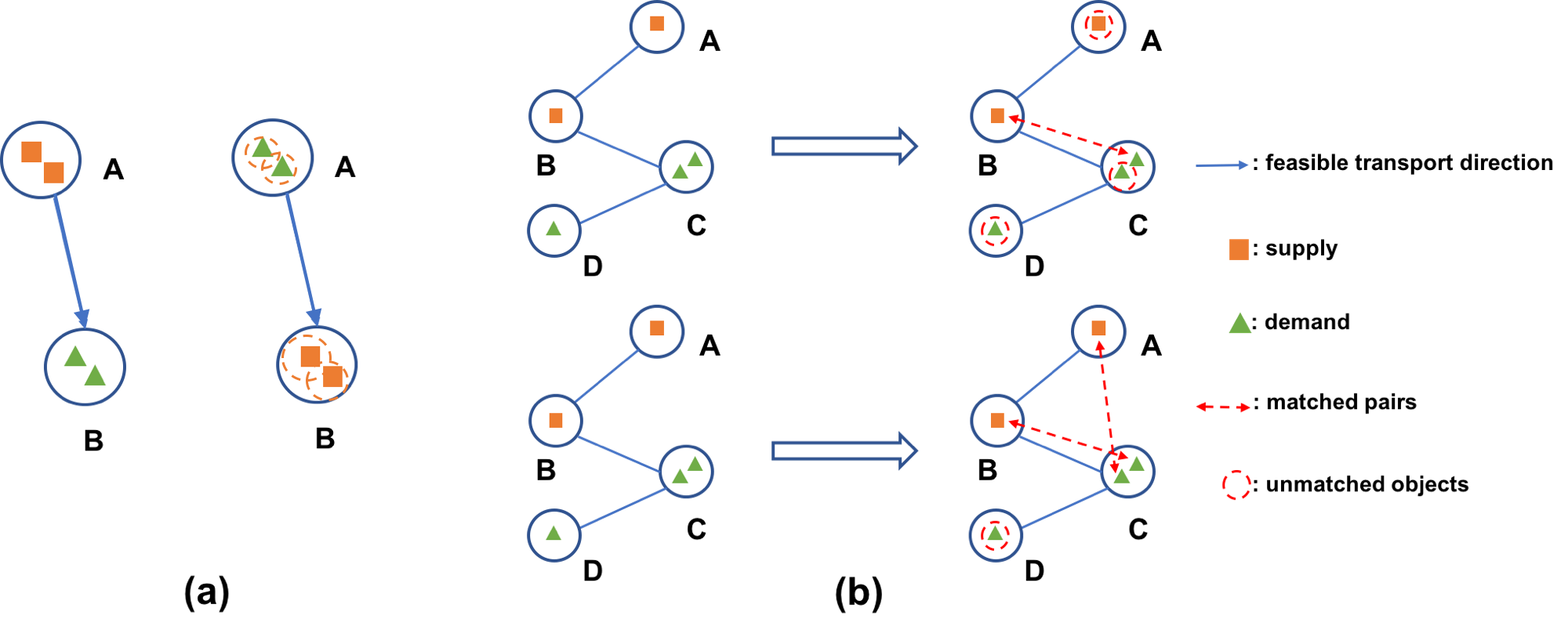} 
\vskip -0.3 in 
	\caption{Examples illustrating the differences between GEM and GWD. Panel (a):   
		For GEM,  the four units can be matched in the left sub-figure,  whereas it is infeasible in the right one.  
		There exits a directed edge from    A to   B, but not from   B to   A. Panel (b): In the top sub-figure, one  `demand' unit at vertex C cannot be matched (the upper line) for GEM since the transport from vertex A to vertex C is not allowed in this case,  whereas in the bottom sub-figure, they can be transported for GWD.  In panel (b), for GEM,  it is assumed that the neighboring set $\mathcal N_i$ only includes the adjacent vertexes of each vertex. 
	}
\vskip -0.2 in  
	\label{compare}
\end{figure}

In our GEM framework, one of the two measures plays the role of 'predator' to move and 'catch' the 'prey', which mimics the general supply-demand system of ride-sourcing platforms. Therefore,    different from the setting of \cite{piccoli2014generalized}, in which both two measures are rescaled, we fix $\nu$ but change $\mu$ only to make the two sides match each other under the asymmetric distance and transport range constraints. 
In Figure \ref{compare}, we  consider two simple examples in order to understand the differences between GEM and GWD.
Moreover,  since we only consider 
the transport from  $\mu$ to $\tilde\mu$ with $\nu$ fixed,   
$\rho_{\lambda} (\mu,\nu | G, C) $ is  generally asymmetric and can be regarded as a restricted GWD.

Besides GEM, 
the optimal solution of $(\tilde\mu, \gamma)$ to (\ref{em}), denoted as $(\tilde\mu_*, \gamma_*)$,  also plays an important role in various two-sided markets,   such as ride-sourcing platforms and E-commerce.  
The $\tilde\mu_*$ 
can be regarded as an optimal dispatch of 
transporting supplies $\mu$ to  match demands $\nu$, whereas $\gamma_*$ is an optimal transport function 
associated with $\rho_{\lambda} (\mu,\nu |G, C)$.  If we vary the area of each vertex 
from the coarsest to the finest scale, then we obtain multilevel GEM and  its transport function. At the finest scale, our GEM reduces to solving an unbalanced assignment problem, so 
$\gamma_*$ is able to capture  the local structure of the optimal transport function. 
In contrast, at a relatively coarse scale,  we obtain a coarse representation of  
the optimal transport function,  reflecting its global patterns.  
We will discuss how to apply GEM to ride-sourcing platforms  in  
Section \ref{sec24}.

Furthermore,  we can simplify $\gamma$ by defining $\gamma = (\gamma_{ij})$ as an $N\times N$ flow matrix with $\gamma_{ij}$ being the transport amount from $v_i$ to $v_j$. Let $\Gamma$ represent  the set consisting of  all the feasible solutions $\gamma$ with all non-negative elements $\gamma_{ij}\geq 0$.
Let $\widetilde{\bm\mu}  = (\widetilde{\mu}_{1},\ldots, \widetilde{\mu}_{N})^T \in R^N$ represent the measure $\mu$ after transporting  $\gamma$ such that $\widetilde{\mu}_{i}=\sum_{v_j \in \mathcal {\mathcal N}_i} \gamma_{ji}$ holds for all $i$.  Thus, our  GEM is equivalent to solving a discrete optimization problem as follows:
\small
\begin{eqnarray}\label{op1}
\rho_{\lambda} (\mu,\nu |G, C) & =& \min_{\gamma \in \Gamma} \{ \parallel \bm{\nu} -  \widetilde{\bm{\mu}} \parallel_1+\lambda \sum_{v_i \in \mathbb V}\sum_{v_j \in \mathbb V} c_{ij}\gamma_{ij}\} \\ 
& \text{subject to } &  \sum_{v_j \in {\mathcal N}_i} \gamma_{ij} = \mu_{i}, \,\,\,\, \sum_{v_j \notin {\mathcal N}_i} \gamma_{ij} = 0, \text{ and}  \,\,\,\, \widetilde{\mu}_{i}=\sum_{v_i \in \mathcal {\mathcal N}_j} \gamma_{ji} \text{ for}\,\,\,\,  \forall v_i \in \mathbb V,  \nonumber 
%& &\text{s.t.} ~~~ { \sum_{v_j \notin {\mathcal N}_i} \gamma_{ij} = 0 \,\,\,\,\,  \forall v_i \in V, \label{Eq1}} \\ 
%& &\text{s.t.} ~~~ \widetilde{d}_{i}=\sum_{v_j \in \mathcal {\mathcal N}_i} \gamma_{ji} \,\,\,\,\,  \forall v_i \in V,   \nonumber 
\end{eqnarray}
\normalsize
where
$\bm{\nu} = (\nu_1,\ldots,\nu_N)^T \in R^N$ and  $\parallel \cdot \parallel_1$ corresponds to  the $L_1$ norm. Moreover,  $\parallel \bm{\nu} -  \widetilde{\bm\mu} \parallel_1$ in (\ref{op1}) is equivalent to the first term of the objective function in (\ref{em}). 
%The global GEM $\rho_{\lambda} (\mu,\nu | G, C)$ includes  the difference between orders in $\mathbf{o}$  and transported drivers in $\mathbf{\widetilde{\mu}}$, the transporting cost $\lambda\sum_{v_i \in V}\sum_{v_j \in V} c_{ij}\gamma_{ij}$, and three sets of constraints.  

There are two key advantages of  using the derived form given in (\ref{op1}) compared to the existing unbalanced optimal transport problem. The first one is that  transport is only allowed between a vertex $v_i$ and its neighboring set $\mathcal N_i$ based on $G$. 
%The second one is that $\lambda$ as a hyper-parameter can balance the transport cost taken to reallocate point masses and the requirement of balancing between two measures $\mu$ and $\nu$. Vertex $v_i$ provides resources to $v_j$ if and only if $\lambda c_{ij}<2$,  which makes a positive contribution to decreasing the objective value. 

The second one is that $\lambda$   can balance the transport cost taken to reallocate point masses and %the requirement of balancing the two measures $\mu$ and $\nu$. 
the requirement of assigning $\mu$ to satisfy $\nu$. The choice of $\lambda$ in practice is data-driven. To ensure that the transport only happens among selected vertex pairs under the optimal transport plan, the theoretical upper bound of $\lambda$ is $2/\max_{v_i\in V, v_j\in \widetilde{\mathcal N}_i} (c_{ij})$,  where $\widetilde{\mathcal N}_i \subset \mathcal N_i$ contains all the neighboring vertexes of $v_i$ that transport from $v_i$.  
%The theoretical upper bound of $\lambda$ is $2/\max_{v_i\in V, v_j\in\mathcal N_i} (c_{ij})$ to make 
%The second one is that the theoretical upper bound of $\lambda$ is $2/\max_{v_i\in V, v_j\in\mathcal N_i} (c_{ij})$ to make GWD meaningful in balancing the transport cost taken to reallocate point masses and the requirement of assigning $\mu$ to satisfy $\nu$ among neighboring vertexes.
%balancing between two measures $\mu$ and $\nu$.
%The second one is that $\lambda$ as a hyper-parameter can balance the transport cost taken to reallocate point masses and the requirement of balancing between two measures $\mu$ and $\nu$. In theory, the upper bound of $\lambda$ is $2/\max_{v_i\in V, v_j\in\mathcal N_i} (c_{ijt})$ to makes GWD meaningful. 
%If we choose $\lambda$ to be  smaller than $2/\max_{v_i\in V, v_j\in\mathcal N_i^1} (c_{ijt})$ for any  $v_i \in V$ and its corresponding $v_j \in \mathcal N_i^1$, 
In this case, the cost of transporting one unit of supply from vertex $v_i$ to $v_j \in \widetilde{\mathcal N}_i$, $\lambda c_{ij}$, is smaller than its contribution to reducing $\parallel \bm{\nu} -  \widetilde{\bm{\mu}} \parallel_1$, which is $2$ ($1$ for $v_i$ and $v_j$, respectively), 
%gain in balancing the demand and supply in 
%the target location $j$ 
when $\nu_{j} - \mu_{j} \geq 1$ and $\mu_{i} - \nu_{i} \geq 1$. Transport from $v_i$ to $v_j$ keeps decreasing the objective value $\rho_{\lambda} (\mu,\nu |G, C)$ until either the balance in the destination vertex $v_j$ or that in the origin vertex $v_i$ is achieved.
%To encourage transporting among all the connected vertexes, we assume that $\lambda$ is smaller than $2/\max_{v_i\in V, v_j\in\mathcal N_i} (c_{ij})$ such that $\lambda c_{ijt}<2$ holds for each $v_i$, $v_j \in \mathcal N_i$ pair.
%throughout the paper. 
In the real world, we usually let $\lambda \max_{v_i\in V, v_j\in \widetilde{\mathcal N}_i} (c_{ij})$ fall into the range $[0.4, 0.5]$ with $c_{ij}$ being the geological distance 
%between $v_i$ and $v_j$ when 
and $\widetilde{\mathcal N}_i$ containing all the first-order adjacent vertexes of $v_i$ in $(G, W, C)$, which can achieve the best performance in some  problems, such as the prediction of order answer rate.

\subsection{Computational Approach} 
\label{sec23}

%We develop an optimization algorithm to solve (\ref{op1}) and compute the optimal solution $\gamma^*$. 
Optimal solution $\gamma^*$ to (\ref{op1}) can be calculated by solving a standard linear programming (LP). 
We will  reformulate (\ref{op1})  as a  LP    problem and then  use a revised simplex method incorporated in a C package GNU Linear Programming Kit (GLPK) to solve   (\ref{op1}).    
We have found that GLPK works pretty well in our real data analyses in Section \ref{sec4}.  

We need to introduce some notations. 
Since the transport range constraints in (\ref{eq8}) impose  $\gamma_{ij}=0$ for $v_j \notin {\mathcal N}_i$, we only need to assign optimal values to $\widetilde \gamma = \mbox{Vec} \{\gamma_{ij}, j \in {\mathcal N}_i\} \in R^{N_0 \times 1}$,  where $N_0 = \sum_{i=1}^N n_i$ 
and  $\mbox{Vec}(\cdot)$   denotes the vectorization of a matrix. With this simplification, the dimension of solvable variables is reduced from $O(N^2)$ to $O(N_0)$, which highly increases the computational efficiency of our algorithm. 
Let $A_1$ and $A_2$ be two $N \times N_0$ matrices. The $i$-th row of $A_1$ consists of $0$'s except the $(\sum_{j=1}^{i-1}n_j +1)$-th to $(\sum_{j=1}^i n_j)$-th elements being $1$. Similarly, all the elements of of $i$-th row of $A_2$ are  zeros except  the $(\sum_{p=1}^{j-1} n_p + q)$-th element being $1$ when grid $v_i$ is indexed by $q$  in  the neighboring set ${\mathcal N}_j$ of vertex $v_j$. 
%Let $\mbox{Vec}(\cdot)$ be  the vectorization of a matrix. 
Let $\widetilde C \in R^{N_0 \times 1}$ be the vector including the unit transport costs for all the corresponding ${\gamma_{ij}}'s \in \widetilde \gamma$. Moreover, we define 
%Let $S \in R^{N \times 1}$ be the $N\times 1$ vector of extra variables introduced for transferring the original $L_1$ norm   problem into a standard LP.  We also define $A$ and ${\bf b}$ as follows: 
\[
A = \begin{bmatrix} 
A_1 & {\bf 0} &  {\bf 0} &  {\bf 0}|| 
A_2 &  {\bf I}_N & -{\bf I}_N& {\bf 0}||
A_2 & - {\bf I}_N &  {\bf 0} & {\bf I}_N
\end{bmatrix}~~~\mbox{and}~~~
{\bf b} = \begin{bmatrix} 
\bm{\mu}^T,
\bm{\nu}^T, 
\bm{\nu}^T
\end{bmatrix}^T,\] 
where $\bm{\mu}=(\mu_1, \ldots, \mu_N)^T$,  $A \in R^{3N \times (N_0 + 3N)}$, ${\bf b} \in R^{3N}$, and ${\bf I}_N$ is an identity matrix. 

The    (\ref{op1}) is equivalent to
$\min
\{ \parallel  \bm{\nu} - A_2 \widetilde \gamma \parallel_1 + \lambda \widetilde C^T \widetilde \gamma\}$  subject to $ A_1 \widetilde \gamma = \bm{\mu}$  and $\widetilde \gamma \geq 0.$ 
Let $S \in R^{N \times 1}$, it can be further  transferred into a standard linear programming (LP)
\begin{eqnarray} 
&& \min\{ \mathbf{1}^TS + \lambda  \widetilde C^T  \widetilde\gamma\}~~~ \text{subject to}  \\
&& A_1 \widetilde \gamma  = \bm{\mu}, ~ A_2 \widetilde \gamma + S \geq \bm{\nu}, ~ A_2 \widetilde \gamma - S \leq \bm{\nu}, ~ \widetilde \gamma \geq 0, ~~
\text{and}~~ S \geq 0.   \nonumber
\end{eqnarray} 
The above LP can be further rewritten as 
\begin{equation}\label{lp}
\min_{X} \{B^T X\}~~~~
\text{subject to}  \,\, A X  ={\bf b}, \,\,\, X \geq 0, 
\end{equation} 
where  $B = (\lambda\widetilde C^T,   {\bf 1}^T, {\bf 0}^T, {\bf 0}^T)^T$  and
$X = (\widetilde \gamma^T, S^T, \bm{w}_{1}^T, \bm{w}_{2}^T)^T$,  in which $\bm{w}_{1}$ and $\bm{w}_{2}$ are vectors of slack variables. 
The  dual of (\ref{lp}) is assigned as 
\begin{equation}\label{dual}
\max_{\textbf{y}\in R^{3N}} \{{\bf b}^T \textbf{y}\}~~~~
\text{subject to}  \,\, A^T \textbf{y}  \leq B, 
\end{equation} 
which further reduces the variable dimension from $N_0 + 3N$ to $3N$.

\subsection{Applications of GEM in Ride-sourcing Platforms}  
\label{sec24}

To calculate GEM, we need to 
build a dynamic  weighted graph structure  over time for each city on the ride-sourcing platform as follows. 
We first  divide  a city  into $|\mathbb V| = N$ non-overlapping hexagons and regard each hexagon as a vertex in $\mathbb V$. Then, 
we set  $\mathcal N_i = \cup_{k=0}^2 {\mathcal N}_i^k$, where ${\mathcal N}^k_i$ includes all the neighboring hexagons  within the $k$-th outer layer of $v_i$ for $k > 1$ and ${\mathcal N}^0_i$ only includes $v_i$ itself.
%and $v_i$ itself. 
A vertex $v_j$ belongs to the $k$-th outer layer of $v_i$ if $k$ steps are required to walk from $v_i$ to $v_j$ on the hexagonal network.  Thus, we determine $G=(\mathbb V, \mathbb E)$.  Second, we set $W_t=(w_{ijt})$, where 
$w_{ijt}$ is  the distance between $v_i$ and $v_j$ in the $t-$th timestamp. Note that $w_{ijt}$  may  vary with  time due to the real-time  locations of drivers and customers.  Third, we compute $C_t=(c_{ijt})$ by using $W_t$  
  in the $t-$th timestamp.  Finally, we obtain the dynamic  weighted graph structure  $(G, W_t, C_t)$. 

We show how to use GEM to address three   operational tasks of interest in ride-sourcing platforms.   
First,  we can  measure the optimal distance between observed dynamic supply and demand networks across time.    
We extract the spatio-temporal data $\textbf{O} =\{ (o_{it})\}_{t}$ and $\textbf{D} = \{(d_{it})\}_{t}$ from the dynamic demand and supply systems, where $o_{it}$ and $d_{it}$ represent demands and supplies at vertex $v_i$  in the $t$-th timestamp, respectively.  
Given $\textbf{O}$ and $\textbf{D}$, we set  $\mu_t=(d_{it})_i$ and $\nu_t=(o_{it})_i$  and use the LP algorithm to calculate 
$\rho(t)=\rho_{\lambda} (\mu_t,\nu_t |G, C_t)$  and its corresponding solution,  denoted as $(\tilde\mu_{t*}=(\tilde d_{it*})_i, \gamma_{t*})$,  in the $t$-th timestamp.

Furthermore,  we introduce an optimal supply-demand ratio at each $v_i$ in the $t-$th timestamp defined as  the ratio of 
$o_{it}$ 
over the 'optimal' supplies $\tilde d_{it*}+ \iota_{\{=\}}(\tilde d_{it*}=0)$, denoted as $
\mbox{DSr}_{it}$, in which we add an extra term $\iota_{\{=\}}(\tilde d_{it*}=0)$ to avoid zero in the denominator. 
Similarly, we can define an optimal supply-demand difference 
as $\mbox{DSd}_{it}=o_{it}-\tilde d_{it*}$ at each $(v_i, t)$.  
It  allows us to create the spatiotemporal map of GEM-related measures $
(\mbox{DSr}_{it}, 
\mbox{DSd}_{it})$. Furthermore, 
we  extend $
(\mbox{DSr}_{it}, 
\mbox{DSd}_{it})$ to a wide timespan $\mathcal T_0$ within a large region $\mathbb V_0 \in \mathbb V$. For instance, we define a  weighted average supply-demand ratio over $\mathbb V_0$ in $\mathcal T_0$  and
a  weighted average absolute  supply-demand difference over $\mathbb V_0$ in $\mathcal T_0$   as follows: 
\begin{equation}\label{ratio_map}
\mbox{DSr}_{\mathcal T_0}(\mathbb V_0) =  \frac{\int_{t \in \mathcal  T_0} \sum_{i \in \mathbb V_0} w_{it} \mbox{DSr}_{it}dt}{\int_{t \in \mathcal  T_0}\sum_{i \in \mathbb V_0} w_{it}dt}~~\mbox{and}~~ 
\mbox{ADSd}_{\mathcal T_0}(\mathbb V_0) =  \frac{\int_{t \in \mathcal  T_0} \sum_{i \in \mathbb V_0} w_{it} |\mbox{DSd}_{it}|dt}{\int_{t \in \mathcal  T_0}\sum_{i \in \mathbb V_0} w_{it}dt}, 
\end{equation}
in which   we set $w_{it}$ as  either $o_{it}$ or $(o_{it}+\tilde d_{it*})/2$ in order to 
highlight vertices with high demands. 
A good market equilibrium in ride-sourcing platforms corresponds to small values of $|\mbox{DSd}_{it}|$  and  
$|\mbox{DSr}_{it}-1|$  across  all $(v_i, t)$. 
Please see    Section \ref{sec4.1}  for details. 

%A fundamental question  of interest for ride-sourcing platforms is 
%\begin{itemize} 
% \item{} how to evaluate whether a ride-sourcing platform achieves  a healthy equilibrium between dynamic supply and demand networks, what we call market equilibrium.  
% \end{itemize} 
% A healthy market equilibrium needs to meet three basic conditions.  The first is that  each  passenger's request  could be easily and quickly satisfied.  Second, the total  idle time for all drivers could be minimal. The last is to minimize various `transporting' costs, such as wait time for riders and pick-up time for drivers, which depend on the size and city geography, weather conditions and market competition. Achieving  such healthy market equilibrium  requires the precise  operation of various platform policies including dispatching, dispositioning, pricing, and subsiding under a given non-platform environment, including weather, holiday, economic level,  and government policy.   
%   how to achieve this equilibrium.  

%The range of $R(t)$ is between $0$ and $1$.  

%Ride-sourcing platforms strive to reduce 

Second, we can use    historical supply-demand information contained in  
$\{(\mbox{DSr}_{it},  \mbox{DSd}_{it}): (v_i, t)\in \mathbb V\times \mathcal T_0\}$ to design   order dispatching policies for large-scale ride-sourcing platforms. 
Order dispatch is an essential component of any ride-sourcing platform for assigning idle drivers to nearby passengers.  
Standard order dispatching approaches    focus on immediate customer satisfaction such as serving the order with the nearest drivers \citep{liao2003real} or the first-come-first-go strategy to serve the order on the top of the waiting list with the first driver becoming available \citep{zhang2016control}. Those greedy methods, however,  fail to account for  the spatial effects of an order and driver (O-D) pair on the other O-D pairs. Thus, they  may not be optimal   from a global perspective. To improve users' experience,  some more advanced techniques strive to balance between small  pick-up distance and 
large  drivers' revenue.   To design better dispatching policy, we will include additional  historical supply-demand network information based on  GEM  to 
delineate its effects   on the average expected gain from serving  current order. Please see Section \ref{sec4.2} for details.

%in the $t$-th timestamp, $\rho(t)=\rho_{\lambda} (\mu_t,\nu_t |G, C_t)$  can be used as a baseline metric for evaluating  the performance of a real matching policy   in any ride-sourcing platform. 

Third,   an important application of  $\{\rho(t)\}$ is to use it as a metric to directly compare two (or more) dispatching policies for ride-sourcing platforms.  
  The key idea is to detect whether there exists a significant difference between two sets of GEMs for two competitive policies under the same platform environment.  Given the joint distribution of demand and supply in the platform, the smaller GEM is, the better many global operational 
metrics, such as order answer rate, order finishing rate, and driver's working time,  are.  Compared with those global operational metrics, GEM is a more direct measurement 
of the operational efficiency for a ride-sourcing platform.    
Please see Section \ref{sec4.3} for details.

\vskip -0.3 in 

\section{Theoretical Properties}  \label{sec3}

In this section, we study the theoretical properties of our  GEM related methods proposed in Section \ref{sec2}, 
most of whose proofs can be found in the supplementary document.

First, we establish the convergence property of   LP (\ref{lp}) for GEM.  

%\noindent To show the primal and dual problems have equal optimal value, we introduce the Strong Duality Theorem.

\begin{theorem}\label{th_3.1} 
	The LP (\ref{lp}) has an optimal basic feasible solution.  Furthermore, 
	if $X$ is feasible for the primal problem (\ref{lp}) and $\textbf{y}$ is feasible for the duality (\ref{dual}), then we have
	\begin{equation}
	\bar{z} = \textbf{y}^T  {\bf b} = \textbf{y}^T A X \leq  B^T  X = z
	\end{equation}
	If either (\ref{lp}) or (\ref{dual}) has a finite optimal value, then so does the other, the optimal values coincide, and  the optimal solutions to both (\ref{lp}) and (\ref{dual}) exist. 
\end{theorem}
%where the following two lemmas are required to complete the proof of Theorem \ref{th_5.1}
%\begin{lemma}[Farkas' lemma] \label{lm2}
%Let $\mathscr{A} \in R^{m\times n}$ and $\mathscr{B} \in R^{m}$ be arbitrary real matrix and vector. Then exactly one of the following two statements is true:
%\begin{enumerate}
%\item There exists an $\mathbf {x} \in R ^{n}$  such that $\mathscr{A} \mathbf {x} =\mathscr{B}$ and $\mathbf {x} \geq 0$
%\item There exists a $\mathbf {y} \in R ^{m}$ such that $\mathscr{A} ^T \mathbf {y} \leq 0$ and $\mathscr{B} ^T \mathbf {y} > 0$
%\end{enumerate}
%\end{lemma}

%\begin{lemma}[\cite{forsgren2008elementary}] \label{lm3}
%Assume that (LP) ($\min \mathscr{C}^T \textbf{x}, \text{ s.t. } \mathscr{A}\textbf{x} = \mathscr{B}, \textbf{x}\geq 0$) is feasible. Then either (LP) has at least one optimal basic solution or there exists a $\textbf{p}$ such that $\textbf{p} \geq 0, \mathscr{A}\textbf{p}=0, \mathscr{C}^T \textbf{p} <0$.
%\end{lemma}

%Proofs of theorems in this section are given in the supplement. The variable dimension of (\ref{lp}), $N_0 + 3N$, increases with $k$, but its duality has fixed dimension $3N$. Solving the duality can highly reduce the variable dimensionality and increase the computation efficiency when $k$ gets extremely large. Below is the dual form of (\ref{lp}):

An implication of Theorem \ref{th_3.1}  is that the LP algorithm for GEM  converges. It demonstrates that there always exit theoretically optimal transport plans (including no transport case) to maximally increase the systematic coherence between the initially unbalanced supplies and demands. However, Theorem \ref{th_3.1} also indicates that the optimal transport plan may not be unique considering the weighted graph structure and initial  supply and demand distributions.

\begin{figure}[h]
	\centering
	\includegraphics[width=0.6\linewidth]{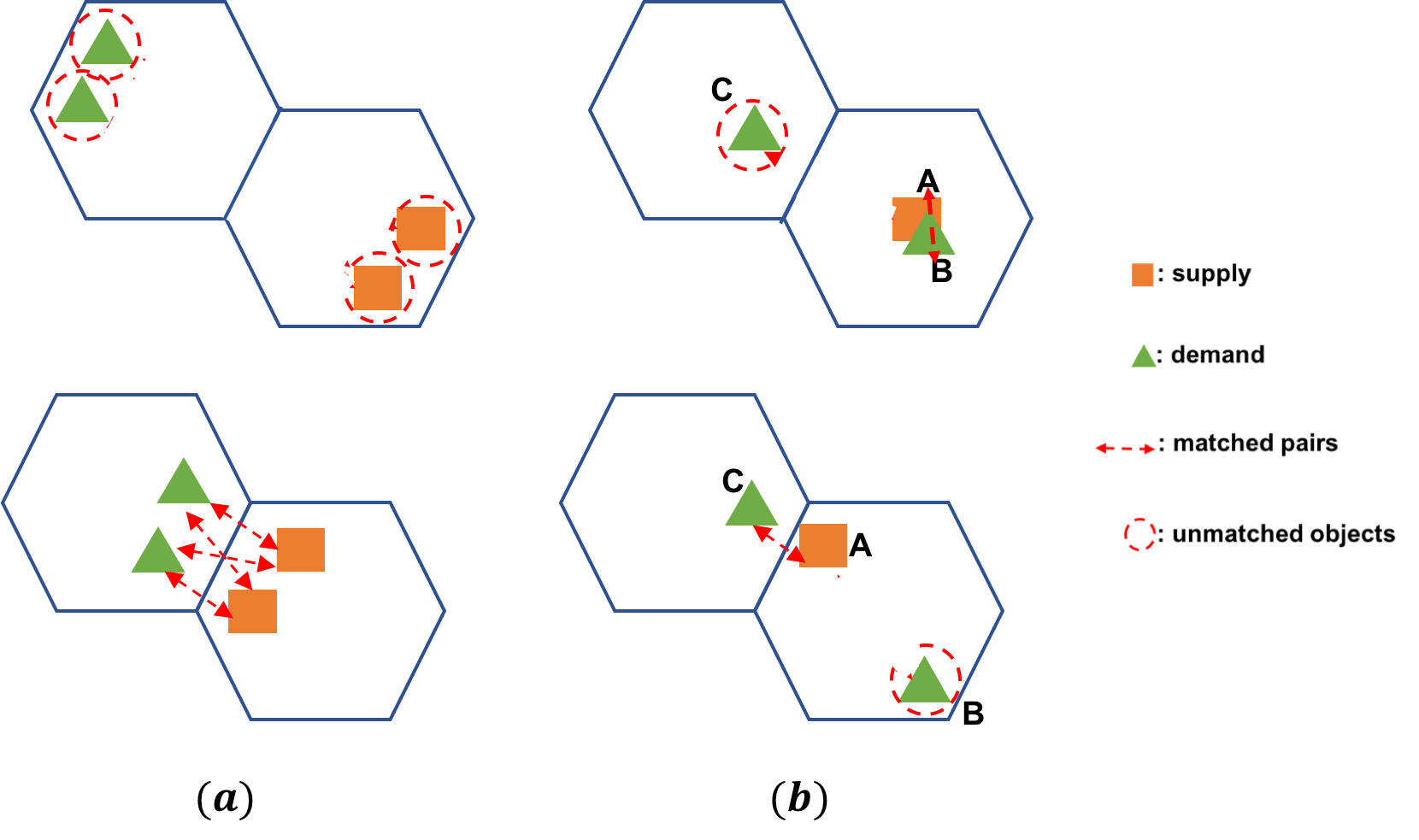}
\vskip -0.2in 
	\caption{Two examples to illustrate the importance of using  random $c_{ij}$.  In Panel (a),  two demands and two drivers  can only be matched in the lower sub-figure since their corresponding pairs of distance  are below a given threshold, whereas it is not the case in the upper sub-figure. 
	%In Panel (b), by letting all $c_{ii}$'s be $0$, A will be matched with B under the solved optimal transport policy to (\ref{op1}). However, pairing A with B may be less costly in practice as the upper sub-figure shows. 
	In Panel (b), supply A is assigned to demand C in the lower sub-figure when the within-grid cost $c_{ii}$ is non-zero, and to demand B when $c_{ii} = 0$, guided by the optimal transport plan with smaller transport costs.
	% instead of A according to the optimal transport with smaller cost when the within-grid cost $c_{ii}$ is non-zero. 
	% in the upper case under the optimal transport, but assigned to demand C in the lower case when the within-grid cost $c_{ii}$ is non-zero. 
	%according to the optimal transport plan to (\ref{op1}) when the within-grid cost $c_{ii}$ is non-zero. 
	%which is optimal in the upper sub-figure other than the lower one. 
	%the assumption that $c_{ii}$ follows a distribution instead of being $0$ allows A to be matched with C from the neighboring vertex if its cost is smaller than that when pairing with B.
	 }
	\label{fig:image3}
\end{figure}

Second, we carry out a probabilistic analysis of our LP (\ref{lp}) for GEM  when $c_{ij}$ follows a distribution. 
Let's start with two motivating examples of ride-sourcing platforms described in Figure \ref{fig:image3}. The $w_{ij}$ represents the geological distance or the traffic time, which may vary between each pair of supply at $v_j$ and demand at $v_i$.  
For the LP defined in (\ref{lp}), it is assumed that each component  of $\widetilde C = (\tilde c_1, \tilde c_2, \ldots, \tilde c_{N_0})^T$ is a non-negative random variable, whereas all elements in $A$ and ${\bf b}$ in (\ref{lp}) are known. Let $z^*$ denote  the minimum value of (\ref{lp}). Since $z^*$ is a function of $\widetilde C$,  it is also a random variable.  We  provide an upper bound for the expectation of $z^*=z^*(\widetilde C)$ below. 

\begin{theorem} [\bf{Expectation Bound}]\label{th_3.2}
	Let $\tilde c_1, \ldots, \tilde c_{N_0}$ be independent non-negative random variables. Suppose there exist $\alpha_1\in (0,\infty)$ and  $\alpha_2 \in (0, 1]$   such that for $l = 1,2,\ldots,N_0$ and all $h>0$ with $P(\lambda \tilde c_l\geq h)>0$,  we have 
	\begin{equation}\label{th6.1.0}
	E(\lambda \tilde c_l|\lambda \tilde c_l\geq h) \geq \alpha_1 \lambda E(\tilde c_l) + \alpha_2 h,  
	\end{equation} 
where the expectation is taken with respect to $\tilde c_l$. 
	{ Let $\{\hat x_1,  \ldots, \hat x_{3N+N_0}\}$
	%$\{\hat x_1,  \ldots, \hat x_{3N+N_0}\}$ 
	be any fixed feasible solution to (\ref{lp}).  
	%$\overline B_l$ denote the expectations of $B_l$ for each $l$, 
%$E(\lambda \tilde c_{[i]})\hat x_{[i]}$ denote the $i$-th largest element of the $E(\lambda \tilde c_j)\hat x_j$'s for $j=1,\ldots,N_0$, and $B_{[i]} = \lambda \tilde c_{[i]}$ be the corresponding $[i]$-th column in $B$, 
We have 
	\begin{equation}\label{th6.1.01}
E(z^*) \leq  {\alpha}_2^{-1}\{\sum_{l=1}^{N_0}(1-\alpha_1\delta_l)E(\lambda\tilde c_l) \hat x_l+\sum_{l=N_0+1}^{N_0+N}\hat x_l\}, 
  %\leq \beta \sum_{i=1}^{N_0 - N} E(B_{[i]})\hat x_{[i]} + \sum_{i = N_0-N+1}^{N_0} E(B_{[i]})\hat x_{[i]} + \sum_{j>N_0} E(B_j)\hat x_j
\end{equation}
where $\delta_l \in [0, 1]$ defined  in the supplementary document is a pre-defined nonnegative constant for each $l \in \{1,\ldots, N_0\}$.}
%If $N_0 \geq 3N$, we have
%	\begin{equation}\label{th_6.1.1}
%	E(z^*) \leq \beta \sum_{i=1}^{N_0 - 3N} E(B_{[i]})\hat x_{[i]} + \sum_{i = N_0-3N+1}^{N_0} E(B_{[i]})\hat x_{[i]} + \sum_{j>N_0} E(B_j)\hat x_j
%	\end{equation}
%{\color{blue} where }

%where $0< \beta < 1$ is a given constant and $B_j$ represents the $j$-th column of $B$.}

%	On the other hand, we assume that $I_{3N}$ contains all the potential basic index sets. For each $r \in I_{3N}$, $B_{(r)}$ and $A_{(r)}$ are the corresponding sub-vector of $B$ and sub-matrix of $A$, respectively. Then we have 
%	\begin{equation}\label{th_6.1.2}
%	\mbox{Var}(z^*) \leq \max_{r\in I_{3N}}\{(A_{(r)}^{-1}b)^T E({B_{(r)}}^T B_{(r)}) A_{(r)}^{-1}b\} - \min_{r\in I_{3N}}[\{E(B_{(r)})A_{(r)}^{-1}b\}^2]    
%	\end{equation}
%	where 
%	\begin{equation}
%	|I_{3N}| = \left\{ \begin{array}{ll}
%	\sum_{x = 3N-N_0}^{2N} \binom{N_0}{3N-x} \binom{3N}{x},   & \mbox{if $N_0 \leq 3N$};\label{th6.1.4}\\
%	\sum_{x=0}^{2N} \binom{N_0}{3N-x}\binom{3N}{x},   & \mbox{otherwise}.\end{array} \right.    
%	\end{equation}
	%\begin{equation}
	%m_0 = \left\{ \begin{array}{ll}
	%N_0 + 3N, & \mbox{if $N_0 < N$};\\
	%4N,  & \mbox{if $ N_0 \geq N$}, \end{array} %\right.
	%\end{equation} 
\end{theorem}
%{\bf Remark}.  

Theorem \ref{th_3.2} has at least two implications.  
First,  condition (\ref{th6.1.0}) holds under some mild conditions.  For instance, it can be shown that  if  $\tilde c_j$ is a  bounded random variable that takes values in   $[c_{j, L},c_{j, U}]$   
such that $P(\tilde c_j>c_{j, L})>0$ and $\lim\inf_{h\rightarrow 0}h^{-1}P(\lambda\tilde c_j<h+\lambda c_{j, L})>0$, then condition (\ref{th6.1.0}) holds. Some   examples of $\tilde c_j$ include uniform, truncated normal, and truncated exponential random variables, among others.  For instance, 
we consider the case that $\tilde c_j$ follows Uniform $[c_{j, L},c_{j, U}]$. It can be shown that  $E(\lambda \tilde c_j|\lambda \tilde c_j\geq h)=0.5(\lambda c_{j, U}+h)$, yielding $\alpha_1=c_{j, U}/(c_{j, U}+c_{j, L})$ and $\alpha_2=0.5$. 
Second,  (\ref{th6.1.01}) gives an upper bound of  the expected value of   $z^*(\widetilde C)$. If we set $\delta_l=0$ for all $l$, then we can obtain a 
larger upper bound  compared with the right-hand side of (\ref{th6.1.01}).  This result generalizes 
an existing result of \cite{Dyer1986} for standard linear programs with random costs under a stronger condition corresponding to $\alpha_1=1$.

%$N_0 \geq 3N$ ($E(z^*) \leq \sum_j B_j \hat x_j$ when $N_0 < 3N$) under condition (\ref{th6.1.0}) since more randomness is included and thus more relaxations can be obtained as $N_0$ grows larger. 
%On the other hand, the results of (\ref{th_6.1.2}) and (\ref{th6.1.4}) consider two cases $N_0 \leq 3N$ and $N_0 > 3N$, respectively. There are more potential combinations of the basic index sets when $N_0 \geq 3N$, which leads to the large size of $I_{3N}$.

Third, we examine the metric properties of $\rho_{\lambda}((\cdot), (\cdot) | G,  C)$ including non-negativity, identity, symmetry, and  the triangle inequality. 

\begin{theorem}\label{th_3.3}
	The operator $\rho_{\lambda}((\cdot), (\cdot) |G,C)$ is a semi-metric such that  it satisfies non-negativity, identity,  and symmetry,  but not necessarily the triangle inequality when 
	(i) $C = (c_{ij}) \in R^{N \times N}$ is symmetric with $c_{ii} =0$ for all $i$; \,\,\,  (ii)  $j \in \mathcal N_i$ if and only if $i \in \mathcal N_j$.  
\end{theorem} 
Theorem \ref{th_3.3} indicates that if $C$ is symmetric, then $\rho_{\lambda} ((\cdot) , (\cdot) |G, C)$ as a semi-metric satisfies three properties including non-negativity, identity, and symmetry. Although the symmetric assumption of $C$ may be incorrect for all vertexes,  it should be valid for most vertexes. Thus, $\rho_{\lambda} ((\cdot), (\cdot)|G, C)$ is approximately a semi-metric. 
	
Fourth, we give the upper and lower bounds of GEM and consider an additivity property in order to better understand how the transport costs and network structures affect GEM.
%{\bf Remark}.  

%Fourth, we examine the XXX  properties of $\rho_{\lambda}((\cdot), (\cdot) | G,  C)$. 

\begin{theorem}\label{thm4.3} The following properties hold:
	
	(i). $\left| |\mu| - |\nu| \right| \leq \rho_{\lambda} (\mu,\nu | G,  C) \leq (|\mu| + |\nu|)$ 
	
	(ii). \textbf{Additivity Property}. For a non-negative $\Delta$,  when $C$ is symmetric, we have 
	%\noindent Given an non-negative integer $\Delta \geq 0$, 
	%let $\mu + \Delta$ be $(\mu_{1}+\Delta, \mu_{2}+\Delta, \ldots, \mu_{N}+\Delta)$ and  $\nu + \Delta$ be $(\nu_{1}+\Delta, \nu_{2}+\Delta, \ldots, \nu_{N}+\Delta)$ as. If $\mathcal C$ is symmetric, then $\rho (\mu, \nu | G, C)$  satisfies the following properties: 
	%\begin{equation}\label{th3eq1}
	\[|\rho_{\lambda} (\mu, \nu +\Delta| G, C) - \rho_{\lambda} (\mu, \nu | G, C)| \leq  N\Delta; \,\,\,\,  |\rho_{\lambda} (\mu +\Delta , \nu | G, C) - \rho_{\lambda} (\mu, \nu| G, C)| \leq  N\Delta. \]
	%and 
	%\begin{equation}\label{th3eq2}
	%-\Delta N \leq \rho (\mu +\Delta , \nu | G, C) - \rho (\mu, \nu| G, C) \leq \Delta N
	%\end{equation}
	
\end{theorem}
%\begin{proof}
%For property (i), we first prove the equality $\left| |\mu| - |\nu| \right| \leq \rho_{\lambda} (\mu,\nu | G, C)$. Without loss of generality, we assume $|\nu| \geq |\mu|$. Also recall that $|\mu| = |\tilde \mu| \leq |\nu|$ by construction. Thus, we have
%\[\rho_{\lambda} (\mu,\nu | G, C) \geq |\nu - \tilde \mu| \geq |\nu| - |\tilde \mu| = |\nu| - |\mu|\]
%To prove the inequality on the right hand side, we choose $\tilde \mu = \nu$ and observe that $\rho_{\lambda} (\mu,\nu | G,  C) = |\nu - \mu| \leq (|\mu| + |\nu|)$ holds. Since $\rho_{\lambda} (\mu,\nu | G,  C)$ is the infimum on all $\mu$ and $\nu$, we have the inequality. The detailed proof of the second property is given in the supplement.
%\end{proof}
Property (i) shows that GEM can be bounded from both above and below. Based on the additivity property, the GEM value can either increase or decrease with one-side node-wise augmentation, which depends on the weighted graph structure and the  distribution of supply and demand. This indicates that applying proper stimulus at selected vertexes is more efficient than globally increasing supply resources.

Fifth, we examine the weak convergence property  of $\rho_{\lambda}((\cdot), (\cdot) | G,  C)$. 

\begin{theorem}\label{thm4.4}\textbf{(Weak Convergence)}
	Let $\{\mu_n\}$ be a sequence of measures on space $\mathbb V$, and $\mu_n, \mu \in M_+(\mathbb V)$. If all the transport costs are bounded, that is $c_{ij} \leq R \,\,$ holds for $\forall  v_i \in \mathbb V$ and $v_j \in \mathcal N_i$, then $\rho_{\lambda} (\mu,\mu_n |G, C) \rightarrow 0$
	%unit cost at all the feasible transport pairs is bounded, i.e. $c_{ij} \leq R, \,\,  \forall \, V_i \in V$ and $V_j \in \mathcal N_i$, we have
	when $\mu_n \rightarrow \mu$ and $\{\mu_n\}$   is tight.  
	%\[\rho_{\lambda} (\mu,\mu_n | G, C) \rightarrow 0 \,\, \text{ is equivalent to } \,\, \mu_n \rightarrow \mu \text{ and } \{\mu_n\} \text{ is tight.}\]
\end{theorem}
Here is an immediate corollary of Theorem \ref{thm4.4}. 
\begin{corollary}
	Let $\{\mu_n\}$ and $\{\nu_n\}$ be two sequences of measures on space $\mathbb V$, and $\mu_n, \nu_n, \mu, \nu \in M_+(\mathbb V)$. If $c_{ij} \leq R$  holds  for $\, \forall \, v_i \in \mathbb V$ and $v_j \in \mathcal N_i$, then we have
	\[\text{if } \mu_n \,\, (resp. \, \nu_n) \rightarrow \mu \,\, (resp.\, \nu) \text{ and } \{\mu_n\}, \{\nu_n\} \text{ are tight, then } \,\rho_{\lambda} (\mu_n,\nu_n |G, C) \rightarrow \rho_{\lambda} (\mu,\nu |G, C).\]
\end{corollary}
Theorem \ref{thm4.4} states that the GEM value goes to $0$ and no transport is required when the initial distributions of $\mu$ and $\nu$ are getting close to each other.

%{\bf Remark}.   

\vskip -0.5 in 

\section{Experiments} \label{sec4}

In this section, we apply  GEM  to the supply-demand diagnostic data set in order to address three  important operational tasks  in ride-sourcing platforms,   including  answer-rate prediction,  the design of order dispatching strategy, and policy assessment. Without special saying, we use the method described in subsection 
\ref{sec24} to construct 
the  dynamic  weighted graph structure  across time  in all these analyses. 
We have released the  `supply-demand diagnostic data set' through the DiDi GAIA Open Data Initiative  at   \url{https://outreach.didichuxing.com/appEn-vue/dataList} and made  the computer codes together with necessary files available at  
\url{https://github.com/BIG-S2/GEM}.   
%A fundamental question of interest   is how to evaluate whether a healthy equilibrium is achieved between dynamic supply and demand systems, what is usually called market equilibrium. A healthy market equilibrium is achieved when passenger's requests could be easily and quickly satisfied, which requires the distance between the two systems to be as minimal as possible.  
%We consider three applications of our GEM, two of which will be presented below and the other is in the supplementary document. 

\subsection{Answer-Rate Prediction}\label{sec4.1}

The data set that we use here includes  both demand and idle driver information from April 21st  to May 20th, 2018 in a large city H. We divide the whole city into $N=800$ non-overlapping hexagonal sub-regions with side length being $1400$m to form the whole vertex set $\mathbb{V}$. { We let the directed edge weight $w_{ij}$ %between two adjacent vertices  
from $v_i$ to $v_j \in {\mathcal N}^1_i$ be the distance between the centers of the two sub-regions, which is $2400$m if $v_j$ can be directly reached by $v_i$ through traffic %there exists any direct road from $v_i$ to $v_j$ 
without first passing through another vertex. Otherwise, $w_{ij} = \infty$.}
%The distance between the centers of two adjacent vertexes is around $2400$m.  
We compute  the numbers of idle drivers and demands in each vertex per  minute and then extract the dynamic supply-demand data set.

\begin{figure}[h]
	\centering
	\includegraphics[width=\linewidth]{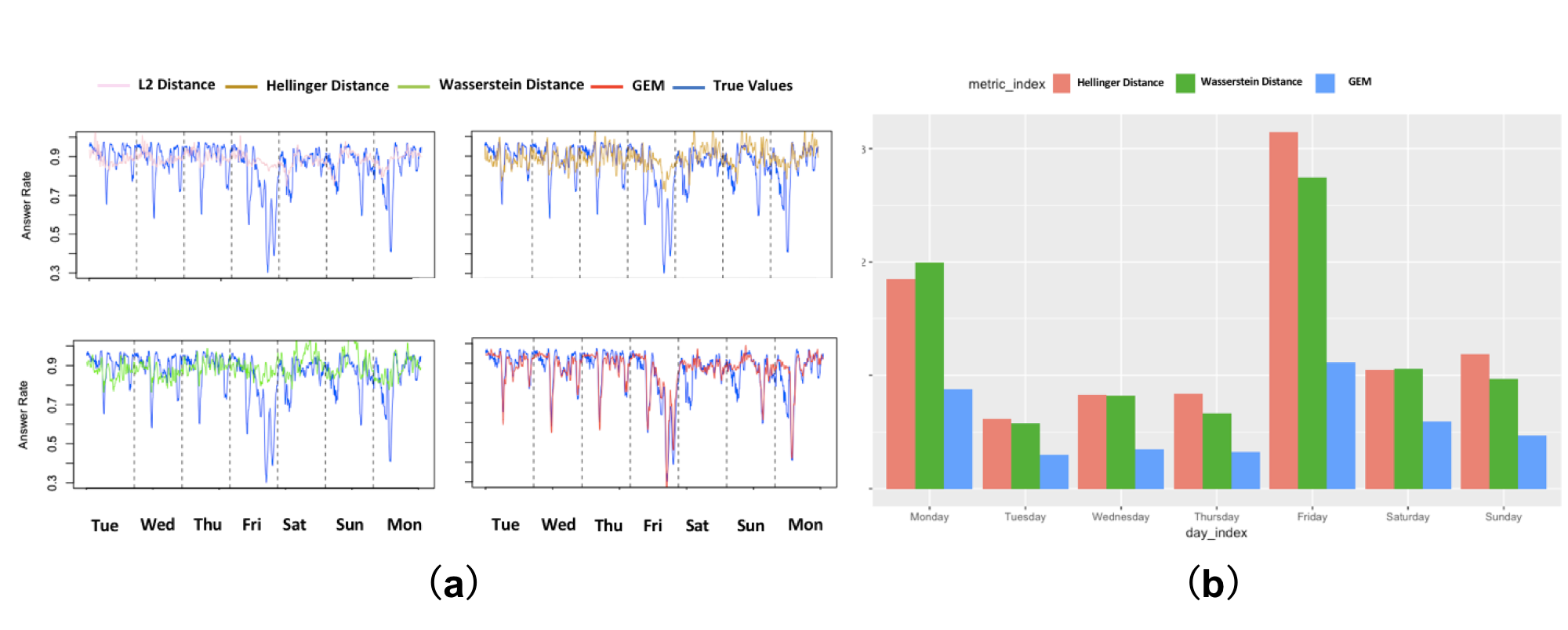}
	\caption{
Results from the answer-rate prediction. Panel (a): comparisons of the log-value of real answer rates obtained   from May 12th to May 18th, 2018
 and their predictive values based on the $L_2$ distance,
the Hellinger distance, the Wasserstein distance, and GEM. Panel (b): comparisons of day-wise RMSEs of answer rate prediction obtained from Monday to Sunday  within the whole city area. }
	\label{fig:prediction}
\end{figure}

The aim of this data analysis is to examine whether the GEM-related measures, such as $\mbox{DSr}_{it}$,  are useful  for predicting  order answer rate  in ride-sourcing platforms. Order answer rate is defined as the number of orders accepted by drivers divided by the total number of orders in a fixed time interval.   
Specifically, we predict  the log-value of order answer rate of the incoming  10 (or 60) minutes by using  historical metric values.  We computed the Hellinger distance, the $L_2$ distance, the 
Wasserstein distance, and  GEM for each 10-minute interval.  
The   $L_2$ distance is calculated by using the numbers of orders and available drivers in all  vertices across 10 consecutive  1-min timestamps.  
The  Hellinger distance is calculated by normalizing the numbers of orders and available drivers in all  vertices and across 10 consecutive  1-min timestamps into probability distributions. For the Wasserstein distance, we first normalize both supplies and demands at each one-minute time interval into two probability distributions and calculate their  corresponding Wasserstein distance. {Subsequently,  we obtain the metric value over each 10-minute interval  by aggregating the Wasserstein distances computed across the 10 included one-minute timestamps by using their corresponding weights $\sum_{v_i \in V} o_{it_k}/\sum_{t_k \in \mathcal T}\sum_{v_i \in \mathbb{V}} o_{it_k}$.}
%every one minute of the 10-minute interval by weights $\sum_{v_i \in V} o_{ik}/\sum_{k \in \mathcal T}\sum_{v_i \in V} o_{ik}$.
%{\bf The Wasserstein distance   is calculated by aggregating the values  calculated at every one minute with weights $\sum_{v_i \in V} o_{ik}/\sum_{k \in \mathcal T}\sum_{v_i \in V} o_{ik}$.} 
For GEM, we compute the supply-demand ratio map  of $\mbox{DSr}_{it}$ per minute  and then we calculate $\mbox{DSr}_{\mathcal T_0}(V)$ for each 10-minute interval.

We split  the supply-demand data set into a training data set consisting  of observations from  April  26th  to May 11th, 2018, and a test data set consisting  of observations from  May 12th  to May 21st, 2018. 
{We use linear regression models to 
predict  the log-value of order answer rate of the incoming $j-$th 10 minutes for $j = 1, \ldots,  6$ by using various historical metric variables of  the previous $p=10$ 10-minute snapshots and those in the same time windows of the previous $5$ days.} 
%as predictors.    

%\begin{figure}
%\centering
%  \includegraphics[width=0.7\linewidth]{plot/ggplot.png}
%  \caption{Answer-rate Prediction:  RMSE by using each metric of different days of the week with the whole city area}
%  \label{fig:RMSE}
%\end{figure} 

We use Mean Absolute Percentage Error (MAPE) and  Root Mean Squared Error (RMSE)  as evaluation metrics to examine the prediction accuracy of all the four compared metrics. 
%evaluate prediction accuracy for all four metrics.  
Table 1 shows their corresponding  RMSE and  MAPE values   based on the test data. Due to the space limitation, we only provide the results corresponding to those at   $t+10$ and $t+60$ minutes, which indicate the  short-term and long-term prediction capacities of all the four metrics. Moreover, we also include the results during the evening  peak hours starting  from 6 pm to 8 pm. For both  the $t+10$ and $t+60$ cases,  GEM significantly outperforms all other three metrics, which may not sufficiently capture the dynamic transport and  systematic  balance of the weighted graph structure.

Figure~\ref{fig:prediction} (a) presents the  real order answer rates and their predictive values in the last 7 test days (Tuesday to Monday)  from May 12th to May 18th   based on  all the four metrics  for the $(t+10)$ case. Compared with  all  other methods,     GEM shows higher consistency between the  true and predicted answer rate values,   especially for some abnormal extreme cases. Furthermore, Figure \ref{fig:prediction} (b) presents  the histograms of RMSEs for  the Hellinger distance, the Wasserstein distance, and 
GEM at each day of the last seven dayes, indicating that GEM outperforms the other two metrics consistently in all seven days. Therefore, 
our GEM is able to capture  the short- and long-term variability within the coherence between the two spatial-temporal systems and has strong  prediction capacity 
for  future answer rates.

\subsection{Order Dispatching Policies}\label{sec4.2}
\begin{figure}
	\centering
	\includegraphics[width=0.65\linewidth]{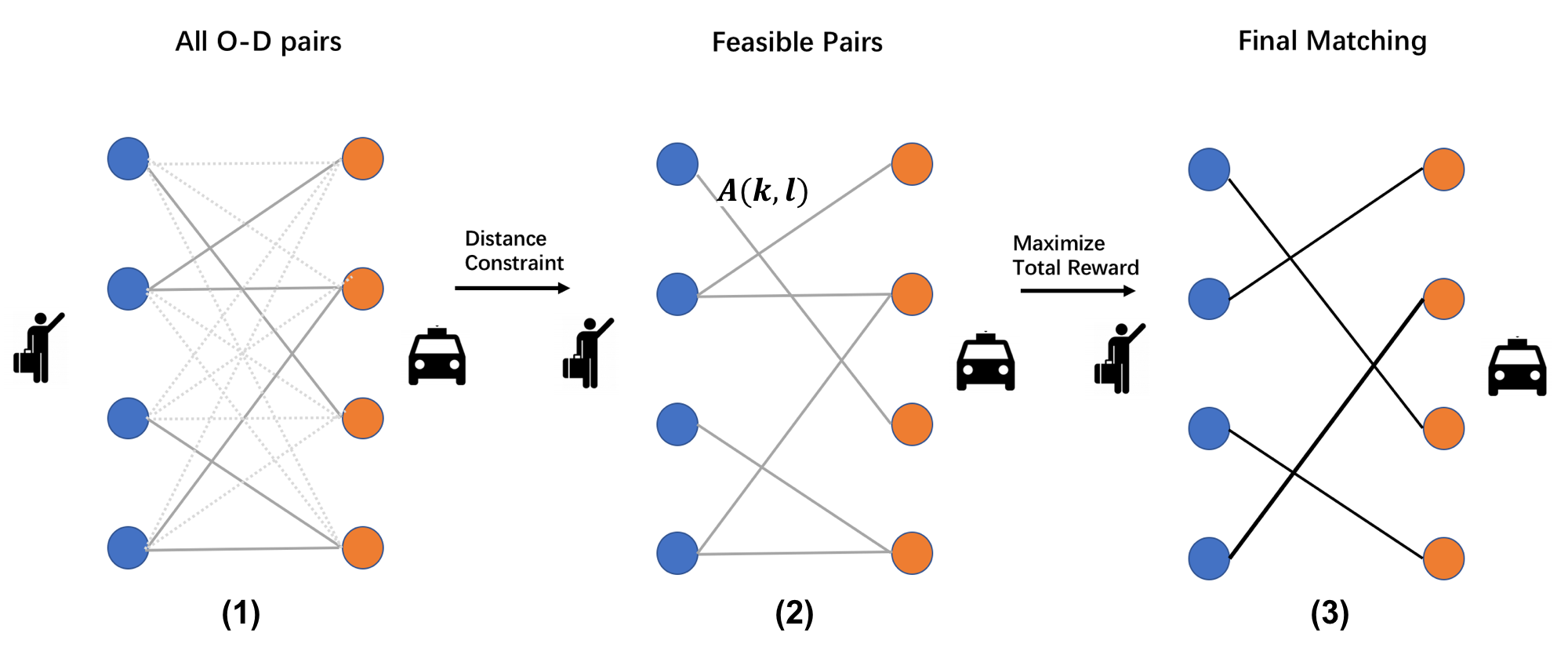}
	\caption{{The order dispatch as a bipartite matching problem: (a) available orders and drivers prepared for pairing;  (b) quantifying all the potential expected earning $A(k,l)$ for all driver-order pairs $(k,l)$ that satisfy the dispatching constraints; and  (c) finding the optimal one-to-one bipartite matching in order  to maximize the total revenue. }}
	\label{fig:dispatch}
\end{figure}

We consider the order dispatching problem of matching $N_o$ orders with $N_d$ available idle drivers, where 
  $N_o$ and $N_d$ denote the total number of orders and that of idle drivers in the current timestamp, respectively.  
The edge weight $A(k,l)$ in the bipartite graph equals to the expected earnings when pairing driver $l$ to order $k$.
%when driver $j$ being assigned to order $i$, which could be driver's revenue, inverse of the pick-up distance or the mixture of the two.  
Let $ x_{kl}$ be    
1  if order $k$ is assigned to driver $l$ and 0 otherwise. 
The global order dispatching algorithm solves a bipartite matching problem as follows:
\begin{eqnarray}
\text{arg}\max_{x_{kl}} \sum_{k=0}^{N_d} \sum_{l=0}^{N_o} A(k,l) x_{kl},  ~  \text{s.t.}  \,\, \sum_{k=0}^{N_d} x_{kl} \leq1 \,\,\,\, \forall l; \,\, \sum_{l=0}^{N_o} x_{kl} \leq1\,\, \forall k;  x_{kl} = 0 \text{ if } c_{kl} > \epsilon~~ \forall k, l. \label{bipartitie} 
\end{eqnarray} 
See Figure~\ref{fig:dispatch} for a graphical  illustration of (\ref{bipartitie}). 
 The constraints ensure that each order can be paired to at most one available driver and similarly each driver can be assigned to at most one order. In practice, only drivers within a certain distance could serve the corresponding orders, which means that $x_{kl}$s' are forced to be $0$ when the distance between order $k$ and driver $l$, denoted as $c_{kl}$,  is beyond the maximal pick-up distance $\epsilon$. The state-of-art algorithm to solve this kind of matching problem is the Kuhn-Munkres (KM) algorithm \citep{munkres1957algorithms}, which will be used to  solve the formulated problem here.

In this paper, we compare three different dispatching policies based on three different formulations of $A(k,l)$. The first one as a baseline only considers the immediate reward of assigning driver $l$ to order $k$, which is defined as
%The simplest form of the edge weight $A(i,j)$ is to only consider the immediate reward of serving driver $j$ to order $i$, such that}
$ 
A^{(1)}(k,l) = \alpha_1 r_k - \alpha_2 c_{kl},
$ 
where $r_k$ is the driver's earning by serving order $k$ and $c_{kl}$ is the pick-up distance between order $k$ and driver $l$. Moreover, $\alpha_1$ and $\alpha_2$ are tuning  parameters such that the two terms are balanced to maximize  drivers' salaries, while reducing  customers'  waiting time.
%minimizing the pick-up distance.
%when the supply is larger than the demand.

{ The second one is given by
$
A^{(2)}(k,l) = \alpha_1 r_k- \alpha_2 c_{kl} + \alpha_3\{\eta^{\Delta t_{lk}} V_1(s'_{lk}) - V_1(s_l)\},     
$ 
where $\eta$ is the discount factor and  an additional term $ \alpha_3\{\eta^{\Delta t_{lk}} V_1(s'_{lk}) - V_1(s_l)\}$ is introduced to enhance the long-term effects of current actions on drivers' future income \citep{xu2018large}. 
%Let $V_1(s_l)$ be the expected earnings from now to the end of the day for driver $l$ located at $s_l = (t, v_{[l]})$, where $t \in T$ is the time index and $v_{[l]} \in \mathbb{V}$ is the current region index of driver $l$.     
%Moreover, 
Let $V_1(s)$ be the expected earnings from now to the end of the day for a driver located at $s = (v, t)$, where $v\in\mathbb{V}$ and $t$ is the current time. Moreover, 
$s_l = (v^{(l)}, t)$ and $s'_{lk} = (v^{(k)}, t+\Delta t_{lk})$ here represent the current spatial-temporal state of driver $l$ and his/her estimated finishing state   when completing serving order $k$, where $v^{(l)} \in \mathbb{V}$ is the current region of driver $l$ before order assignment and  $v^{(k)} \in \mathbb{V}$ is the destination region of order $k$ and $\Delta t_{lk}$ denotes the total time required for driver $l$ to finish the whole process of serving order $k$. 
%where $v^{(k)}$ is the destination region that order $k$ goes to and $\Delta t_{lk}$ denotes the total time required to finish the whole serving process for order $k$.   
If a driver becomes available to a new order immediately after finishing the ongoing one, then  $\eta^{\Delta t_{lk}} V_1(s'_{lk}) - V_1(s_l)$ is the extra future earning for driver $l$ by serving order $k$ other than staying idle.

The third one is given by 
$
A^{(3)}(k,l) = A^{(2)}(k,l)+ \alpha_4  \{\eta^{\Delta t_{lk}}V_2(s'_{lk}) - V_2(s_l)\}, 
$
%- v_1(s))  \nonumber \\ 
%&+ &\alpha_4 (v_{2t}(s'_{ij}) - v_{2t}(s_j)),   \label{Eq3}
%\end{eqnarray} 
%where $v_{2t}(s_j)= o_{g\mathcal T_0} - \tilde{d}^*_{g\mathcal T_0}$ 
%represents the supply-demand ratio 
%in a time interval $\mathcal T_0$ at grid $v_g$. Higher $V_2$ value is equivalent to lower expected demand satisfactory rates in the target spatial-temporal point. 
where    $ \alpha_4  \{\eta^{\Delta t_{lk}}V_2(s'_{lk}) - V_2(s_l)\}$ is   introduced to  balance the supply-demand coherence. 
Moreover,  $V_2(s) = \nu_{t}(v) - \tilde \mu_{t}(v)$  at $s=(v, t)$ is calculated from  GEM in  (\ref{op1}). } Using $V_2(\cdot)$ increases the probability that customers' requests can be quickly answered by nearby drivers, whereas
 $V_1(\cdot)$ ignores the interaction effects when multiple drivers are heading to the same location.
%making the same decision. 
%In this case,  the supply-demand relationship in the destination grid may influence the marginal value of an extra assignment. 
Thus,  when the future demand has already been fulfilled by drivers re-allocated by previous completed servings, assigning more drivers  might decrease   $V_1(\cdot)$ in the target location.

We use a comprehensive and realistic dispatch simulator designed for recovering the real online ride-sourcing system to evaluate the three dispatching policies. 
The simulator models the transition dynamics of the supply and demand systems to mimic the real on-demand ride-hailing platform. The order demand distribution of the simulator is generated based on historical data.  The driver supply distribution is initialized by historical data at the beginning of the day, and then evolves following the simulator's transition dynamics (including drivers getting online/offline, driver movement with passengers and idle driver random movement) as well as the order dispatching policies. The differences between the simulated results and the real-world situation is less than 2\% in terms of some important metrics, such as drivers' revenue, answer rate,  and idle driver rate.

To compare the three dispatching policies, 
we randomly  selected a specific  city S, which usualy has in total 150, 000 to 200, 000 ride demands per day. We still divide the whole city area into $N= 800$ hexagonal vertices and use the geological distance between two nearby grids to be the edge weights. Furthermore, three different days including  2018/05/15 (Tuesday),  2018/05/18 (Friday), and 2018/05/19 (Saturday) were analyzed since the global order answer rates on weekday are usually much lower than those at weekend by looking at the historical data. Both $V_1(\cdot)$ and $V_2(\cdot)$ values were obtained by taking the average of the same weekday or weekend from the previous four weeks since the platform has significant weekly periodicity. The length of time intervals that we used to compute $V_1(\cdot)$ and $V_2(\cdot)$ was set to be 
$\mathcal T_0 = 10$ minutes so that all the action windows inside share the same $V_1(\cdot)$ and $V_2(\cdot)$ values. Specifically, $V_2(\cdot)$ is achieved by aggregating the $|\mathcal T_0|$ continuous $(\nu_{kt} - \tilde \mu_{kt})$s'.  
	%The ratio value $m_{t,g}$ is achieved by taking the average of the same weekday or weekend from the previous four weeks since we find that the matching degree of demand and supply has weekly pattern when the order dispatching policy is consistent by analyzing the historical data. The length of the timestamp we use for the value function $v_{1t}$ \citep{xu2018large} and our ratio value $v_{2t}$ is $\mathcal T_0 = 10$ min. Our ratio measurement is computed every one minute and aggregated within each 10 minutes as Section \ref{sec9.3} does. 
	We applied the three dispatching policies with different edge weights to the simulator even based on  the same initial input and transition dynamics. We set $\alpha_1=1$ and $\alpha_2=0.001$ to rescale the order price $r_k$ and the pick-up distance $c_{kl}$ into comparable ranges.  The $r_k$ contributes more to the variations of $A(k, l)$ because of the constrained pick-up distance ($c_{kl} \leq \epsilon$).
Furthermore, we perform grid search for a wide range of $(\alpha_3; \alpha_4)$ combinations to find its optimal solution,  denoted as 
$(\alpha_3^*, \alpha_4^*)$, that maximizes average drivers’ revenues  for weekdays and weekends in the simulator. 
%Specifically, we fix the weight parameters in (\ref{Eq2}) and (\ref{Eq3}) to be $\alpha_1=1$ and $\alpha_2=0.001$. 
Specifically, we fixed $\alpha_4 = 0$ first and  use the bisection method to obtain a rough value range of length 0.1 for $\alpha_3$ with its initial start being $[0,1]$. 
 Then we apply the grid search method to increase 0.01 amount for $\alpha_3$ each time within the value range until finding the optimal $\alpha_3^*$ corresponding to the largest averaged drivers’ revenue. Subsequently, we fix $\alpha_3^*$ and do the similar grid search to get the optimal $\alpha_4^*$.

Tables~\ref{tbl:order_dispatch} and \ref{tbl:order_dispatch_2} summarize the collected results corresponding to the baseline policy, $A^{(2)}(k,l)$ with the optimal $\alpha_3$, and our approach with different $\alpha_4$ values. It reveals that the order dispatching policy  based on $A^{(3)}(k,l)$ could achieve higher drivers' revenue and answer rate compared with the other two policies. The optimal $\alpha_3$ is achieved at $0.54$, $0.61$, and $0.52$ for 2018/05/15, 2018/05/18, and 2018/05/19, respectively. In 2018/05/15 and 2018/05/19, we obtain  a smaller optimal $\alpha_3$ since a higher coherence between supplies and demands 
is achieved under the baseline policy ($\alpha_3 = \alpha_4 = 0$) than that of 2018/05/18, which indicates that the supply-demand relationship is more related to the policy efficiency than the weekday/weekend status. Moreover, the supply abundance in 2018/05/15 and 2018/05/18 results in a higher order answer rate but a smaller optimal $\alpha_4$. Compared to the policy corresponding to $A^{(2)}(k,l)$, adding the GEM-related measurements increases the expected whole-day answer rate and drivers' revenue in more than $1\%$. It may  indicate  that the supply-demand difference may affect the expected future gain of a marginal driver.

In practice, we first 
 perform grid search for a wide range of $(\alpha_3, \alpha_4)$ combinations to find its optimal solution  $(\alpha_3^*, \alpha_4^*)$  that maximizes average drivers’ revenues for some representative days in the simulator. Then we fine-tune the parameters via on-line A/B testing, and apply the policy in the real-life dispatching system. The value functions $V_1(\cdot)$ and $V_2(\cdot)$ are updated when the new policy being employed for a period of time, and $\alpha_3^*$ and $\alpha_4^*$ are re-tuned in the real environment to achieve the optimal efficiency. 

\subsection{Policy Evaluation}\label{sec4.3}

\begin{figure}[h]
	\centering
	\includegraphics[width=1.05\linewidth]{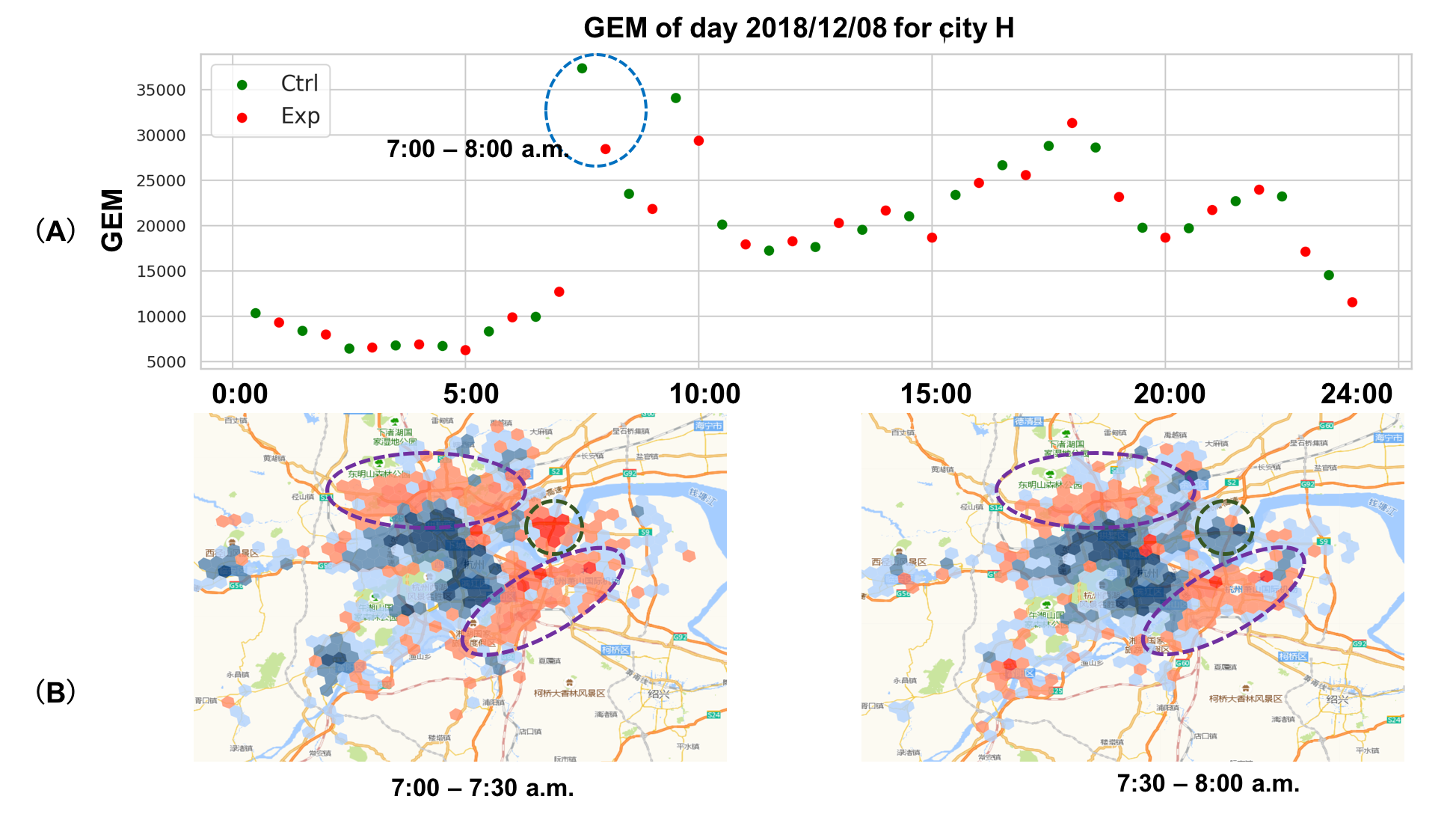}
	\caption{Results from the policy evaluation. (A) The GEM values  of a randomly selected day on 2018/12/08 for city H at 30-min scale.  Green and red points represent the GEM values generated by the baseline (control) and new (experimental) policies, respectively. In particular, we mark the time period 7:00 to 8:00 a.m. by a blue circle, which demonstrates a significant reduction of GEM value when changing the policy from the the control one to the experimental one.   (B) The  heatmaps of vertex-wise supply-demand difference $\mbox{DSd}_{it}$ of city H under the control    and experimental policies  within the 30-min time window from 7:00 to 7:30 am   
and 
 that from 7:30 to 8:00 am,  respectively.   Hexagons in red and blue colors represent the locations with positive and negative $\mbox{DSd}_{it}$, respectively,  and  a deeper color corresponds to a big 
 $|\mbox{DSd}_{it}|$ value.  
}
	\label{gem_compare}
\end{figure} 

% If our proposed GEM can sufficiently represent the coherence between the demands and supplies on the platform, it should be highly correlated with some important indexes such as order answer rates, order finishing rate and driver working time as Section \ref{sec4.1} illustrates. Therefore, we may conclude that A policy outperforms B policy if the EM values significantly decrease.

We conduct an experiment using another supply-demand data set of the same city H  from December 3rd to December 16th, 2018  in order to compare the effectiveness between two order dispatching policies. We executed them alternatively on successive half-hourly time intervals. Moreover, we start with the baseline policy being in the first half hour and change the policy every half hour through the whole day and reverse their order in another day.  We include an A/A test, which compares the baseline policy against itself, by using the historical data obtained from November 12th to November 25th as a direct comparison. 
We calculate GEM within each time window of 30 minutes as follows.  
There are  in total $M_T = 48$ time intervals  per day. To obtain GEM in each time interval $\mathcal T$,  we 
 aggregate  30 GEM values, each of which is calculated within the 1-min timestamp,  by using  normalization weights $o_{it}/\sum_{t\in \mathcal T}\sum_{v_i \in \mathbb{V}} o_{it}$.

We first need to introduce some notation. 
 We denote $y_m(t_k)$ as the aggregated GEM value  and  use $x_{m}(t_k)$ to denote a $2\times 1$ vector of predictors, which   are not strongly influenced  by  order dispatching policy, including the total number of demands and the  total supply time of all drivers in the $k$-th time interval of day $m$ for $k = 1,  \ldots, M_T$ and 
$m=1, \ldots, M_D$. Let   $a_m(t_k) = 1$ if the new policy is used and $=-1$ otherwise. { To examine the marginal effect of policy  on GEM, we consider the following regression model:
\begin{equation}\label{GEEtr}
%y_i(t_k) = \alpha(t_k) + \beta(t_k)^T\{x_i(t_k) - \overline x(t_k)\} + \gamma(t_k) a_i(t_k) + e_i(t_k),
\begin{aligned}
y_m(t_k) &= \beta_0(t_k) + \beta_1(t_k)^T\{x_m(t_k) - \overline x(t_k)\} + \beta_2(t_k) a_m(t_k) + \eta_m(t_k) + \varepsilon_m(t_k), 
\end{aligned}
\end{equation}
where $\beta(t_k) = (\beta_0(t_k), \beta_1(t_k)^T, \beta_2(t_k))^T$ is a vector of regression coefficients at $t_k$,  and $\overline x(t_k)$ is  the sample mean of all $x_i(t_k)$s for $k = 1,  \ldots, M_T$. In addition, we assume that $\eta_m = (\eta_m(t_1), \ldots,\eta_m(t_{M_T}))^T$ and $\varepsilon_m = (\varepsilon_m(t_1),\ldots, \varepsilon_m(t_{M_T}))^T$ are $M_T\times 1$ vectors of random errors, following mutually independent multivariate Gaussian distributions $N(\textbf{0}, \Sigma_{\eta})$ and $N(\textbf{0}, \sigma^2_\varepsilon \cdot {\bf I}_{M_T})$, where $\Sigma_{\eta}$ is an $M_T\times M_T$ matrix and $\sigma^2_\varepsilon$ is a positive scalar.  
%are mutually independent and identical copies of $\textbf{SP}\left\{0, \Sigma_{\eta}(t, t^{\prime})\right\}$ and $\textbf{SP}\left\{0, \sigma_{\varepsilon}^2(t)I(t^{\prime}=t)\right\}$
%respectively, where $\textbf{SP}\left\{\mu, \Sigma\right\}$ denotes a stochastic process with mean function $\mu(t)$ and covariance function $\Sigma(t, t^{\prime})$, and $I(\cdot)$ is the indicator function of an event.
%$(\alpha(t_k), \beta(t_k)^T, \gamma(t_k))^T$       
%is a vector of regression coefficients at  
%$t_k$. 
%Moreover,   $e_i = (e_i(t_1),\ldots,e_i(t_{M_T}))^T$ is a $M_T\times 1$ vector of random errors 
%and it is assumed that $e_i$ follows a multivariate Gaussian distribution $N(\textbf{0}, \Sigma)$, where $\Sigma$ is an $M_T\times M_T$ matrix. 
We are interested in testing the following null and alternative hypotheses: 
\begin{equation}\label{hypothesis}
H_0: \int_{0}^{M_T}\beta_2(t)dt = 0~~~\mbox{v.s.} ~~H_1: \int_{0}^{M_T}\beta_2(t)dt\not=0,  
\end{equation}
where $\int_{0}^{M_T}\beta_2(t)dt\approx \sum_{k=1}^{M_T} \beta_2(t_k)\Delta t_0$ denotes the average treatment effect per day, in which $\Delta t_0$ is the length of each time interval.  
%We propose a novel estimation procedure by using Generalized Estimating Equations (GEE) \citep{liang1986longitudinal, mancl2001covariance} for the above testing problem (\ref{hypothesis}). 
We propose a joint estimation procedure based on Generalized Estimating Equations (GEE) to iteratively estimate  all unknown parameters until a specific convergence criterion being reached \citep{liang1986longitudinal}. Subsequently, we compute the $t-$test statistic associated 
with  the average treatment effect per day and its 
corresponding  one-sided (or two-sided) $p-$value  \citep{mancl2001covariance}.  

} 
%use Generalized Estimating Equations (GEE) \citep{liang1986longitudinal, mancl2001covariance} to estimate all the related parameters in the above regression model.}  

Furthermore, we consider three global operational metrics including the order answer rates, order finishing rate, and gross merchandise value (GMV) as $y_m(t_k)$ in model (\ref{GEEtr}).  We fit the corresponding  three regression models in order to study whether the new dispatching policy significantly improves the ride-sourcing platform at the operational level.

 Table \ref{abtest} summarizes all regression analysis results  for both the  A/A and A/B experimental designs.  
We can see that in the A/B experimental design,  there exists a significant increase in the mean answer rate, finishing rate and gross merchandise value when replacing the old policy by the new one since all the $p-$values associated with the average treatment effect   are smaller than $10^{-3}$. The new policy can also significantly reduce the GEM value ($p-$value smaller than $0.05$), which agrees with our assumption that GEM can sufficiently quantify the supply-demand relationship and subsequently affect the examined platform indexes.  
In contrast,   Table \ref{abtest} shows that in the A/A experimental design, all the four metrics do not show significant treatment effect at the significance level of $5\%$. 

Figure \ref{gem_compare}(A) presents the GEM value at in total 48 30-min time windows on December 3rd, 2018  for the  A/B experimental design.   
We observe a significant reduction of GEM value when changing the policy from the the control one to the experimental one during the time period  from 7:00 to 8:00 a.m.  Figure \ref{gem_compare}(B) presents the heat maps of vertex-wise $\mbox{DSd}_{it}$ within the same time period
under the control and experimental policies, respectively. 
 The customer requests in three selected regions marked by green and purple circles were satisfied by the  drivers in nearby regions,  resulting  in the higher supply-demand coherence and thus a smaller GEM value.

%\begin{figure}
%	\centering
%	\includegraphics[width=1\linewidth]{plot/ratio.png}
%	\caption{Comparison of real order answer rate with $R(t)$}
%	\label{fig:ratio}
%\end{figure}

\bibliographystyle{agsm}
\bibliography{YingPaper2015,ZhuShenR01July11_4}
\newpage

\begin{table*}[htbp]
	\centering 
	\caption{Results from the answer-rate prediction. Comparisons of  Hellinger, $L_2$-distance, Wasserstein, and GEM in predicting answer rate 
at   $t+10$ and $t+60$ minutes. Peak hour  denotes the time   from 6 pm to 8 pm.   MAPE and RMSE denote  the  mean absolute percentage error  and  root mean squared error, respectively.   }
	\begin{tabular}{cccccccc}
		\hline\hline 
		&    &   &  Hellinger   & L2-distance & Wasserstein & GEM \\\hline
		t+10  & All time & RMSE  & 0.1362   & 0.1496   &  0.1273 & \bf{0.0552} \\
		&          & MAPE  &  0.0801  & 0.0891   &  0.0718  & \bf{0.0338}\\
		& Peak hour& RMSE  & 0.2219  & 0.2187   &  0.2088  &\bf{0.0614}  \\ 
		&     & MAPE  &  0.1494  & 0.1457   & 0.1089 &  \bf{0.0422} \\\hline
		t+60  &  All time & RMSE  & 0.1522   & 0.1552   &  0.1413 & \bf{0.1130}\\
		&           & MAPE  &  0.0828  & 0.0868   & 0.0859  & \bf{0.0620} \\
		& Peak hour & RMSE   & 0.2395  & 0.2565  &  0.2222 & \bf{0.1530}  \\
		&       & MAPE       &  0.1077  & 0.1159  &  0.1317 &  \bf{0.0728} \\
		\hline\hline  
	\end{tabular}
	\label{tresult}
\end{table*} 

\begin{table}
	\centering
	\caption{Results from order dispatching policies. Comparisons of the three policies  with respect to two evaluation metrics including the drivers' revenue and the global answer rate using the simulator for city S on two selected Weekdays. 
	%two selected days. 
	The rows with $(\alpha_3, \alpha_4)=(0, 0)$ correspond to 
the first (or baseline) policy, those with $\alpha_4=0$ and $\alpha_3\not=0$ correspond to the second policy, and 
all other rows correspond to the third policy. { The numbers in the parentheses denote the relative improvement of  the corresponding policy over the baseline policy for each evaluation metric.  } }\label{tbl:order_dispatch}
	\begin{tabular}{cccc}
		\hline\hline 
		$\alpha_3$& $\alpha_4$ & Drivers' Revenue (Yuan)  & Order Answer Rate\\
		\hline
		\multicolumn{4}{c}{
			2018/05/15 (Tuesday)} 
		\\
		\hline 
		0       &  0  &1191316 &  0.737\\
		0.54  &  0  &  1227175(+3.01\%)  &  0.760(+3.12\%) \\
		0.54  &  6  &  1235037(+3.67\%)  &  0.761(+3.28\%)\\
		0.54  &  7  &  1236824(+3.82\%)  &  0.763(+3.54\%)\\
		0.54  &  8  &  \textbf{1240518(+4.13\%)}  &  0.765(+3.82\%)\\
		0.54  &  9  &  1238850(+3.99\%)  &  0.764(+3.66\%) \\
		0.54  & 10 &  1231702(+3.39\%)  &  0.761(+3.26\%) \\
		\hline
		\multicolumn{4}{c}{
			2018/05/18 (Friday)} 
		\\
		\hline 
		0       &  0  &13943666 &  0.539\\
		0.61  &  0  &  14701230(+5.43\%)  &  0.557(+3.34\%) \\
		0.61  &  6  &  14845486(+6.47\%)  &  0.561(+4.08\%)\\
		0.61  &  7  &  14858400(+6.56\%)  &  0.560(+3.90\%)\\
		0.61  &  8  &  14865454(+6.61\%)  &  0.560(+3.90\%)\\
		0.61  &  9  &  \textbf{14867573(+6.63\%)}  &  0.560(+3.90\%) \\
		0.61  & 10 &  14823121(+6.31\%)  &  0.557(+3.34\%) \\
		\hline
		%\multicolumn{4}{c}{
		%	2018/05/19 (weekend)} 
		%\\
		%\hline 
		%0     & 0   &  13507568 &  0.745\\
		%0.52  &  0  &  13886185(+2.80\%)  &  0.768(+3.09\%) \\
		%0.52  &  6  &  14034453(+3.90\%)  &  0.774(+3.89\%)\\
		%0.52  &  7  &  14008847(+3.71\%)  &  0.772(+3.62\%)\\
		%0.52  &  8  &  \textbf{14043995(+3.97\%)}  &  0.773(+3.76\%) \\
		%0.52  &  9  &  13996021(+3.62\%)  &  0.770(+3.36\%) \\
		%0.52  & 10 &  13934895(+3.16\%)  &  0.768(+3.09\%) \\
		%\hline
		\hline 
	\end{tabular}
\end{table}

\begin{table}
	\centering
	\caption{Results from order dispatching policies. Comparisons of the three policies  with respect to two evaluation metrics including the drivers' revenue and the global answer rate using the simulator for city S on a selected Weekend. }\label{tbl:order_dispatch_2}
\begin{tabular}{cccc}
		\hline\hline 
		$\alpha_3$& $\alpha_4$ & Drivers' Revenue (Yuan)  & Order Answer Rate\\
		\hline
		\multicolumn{4}{c}{
			2018/05/19 (Saturday)} 
		\\
		\hline 
		0     & 0   &  13507568 &  0.745\\
		0.52  &  0  &  13886185(+2.80\%)  &  0.768(+3.09\%) \\
		0.52  &  6  &  14034453(+3.90\%)  &  0.774(+3.89\%)\\
		0.52  &  7  &  14008847(+3.71\%)  &  0.772(+3.62\%)\\
		0.52  &  8  &  \textbf{14043995(+3.97\%)}  &  0.773(+3.76\%) \\
		0.52  &  9  &  13996021(+3.62\%)  &  0.770(+3.36\%) \\
		0.52  & 10 &  13934895(+3.16\%)  &  0.768(+3.09\%) \\
		\hline\hline 
	\end{tabular}
\end{table}

\begin{table*}[htbp]
	\centering
	\caption{Results from the policy evaluation: relative improvement and two-sided $p-$value of average treatment effects  for the  A/A and A/B experiments}
	\begin{tabular}{cccc}
		\hline\hline 
		Experiment Design  & $y_m(t)$ & Relative Improvement($\%$) & $p-$value \\\hline
		& Answer Rate & 0.76 & 1.16e-12 \\
		A/B & Finish Rate & 0.36 & 4.32e-3 \\
		& GMV & 0.86 & 2.91e-6 \\
		& GEM & -0.80 & 4.06e-2 \\\hline
		& Answer Rate & 0.01 & 0.96 \\
		A/A & Finishing Rate & 0.01 & 0.96 \\
		& GMV & -0.08 & 0.72 \\
		& GEM & -0.25 & 0.43\\\hline
		\hline\hline  
	\end{tabular}
	\label{abtest}
\end{table*} 

\newpage 

%\bigskip
%\begin{center}
%{\large\bf SUPPLEMENTARY MATERIAL}
%\end{center}

%\begin{description}

%\item[Title:] Brief description. (file type)

%\item[R-package for  MYNEW routine:] R-package ?MYNEW? containing code to perform the diagnostic methods described in the article. The package also contains all datasets used as examples in the article. (GNU zipped tar file)

%\item[HIV data set:] Data set used in the illustration of MYNEW method in Section~ 3.2. (.txt file)

%\end{description}

\end{document}